\documentclass[11pt]{amsart}
\usepackage{graphicx, amssymb}

\newtheorem{defn}{Definition}
\newtheorem{thm}{Theorem}
\newtheorem{lemma}{Lemma}
\newtheorem{prop}{Proposition}
\newtheorem{ex}{Example}
\newtheorem{cor}{Corollary}

\begin{document}

\title{A cut-and-paste approach to contact topology}

\author{William H. Kazez}

\address{University of Georgia, Athens, GA 30602}

\email{will@math.uga.edu}

\urladdr{http://www.math.uga.edu/\char126 will}



\thanks{Supported in part by NSF grant DMS-0073029.}

\date{September 26, 2002.}

\maketitle

\vspace{-.2truein}

\begin{center}
{\includegraphics[height=.8cm]{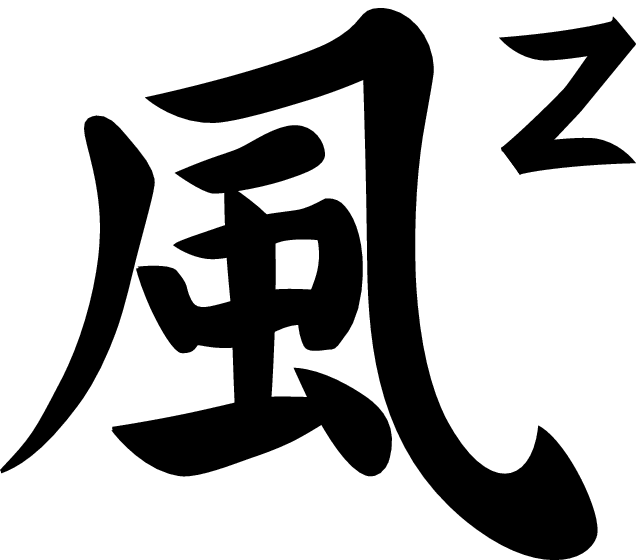}}
\end{center}

\tableofcontents

\section{Convex Surfaces}

Unless otherwise stated,  $M$ will denote a compact  orientable 3-dimensional manifold and may have nonempty boundary.

\begin{defn}
A (positive) {\rm contact structure}, $\xi$ on $M$, is a smooth 2-plane bundle
$\xi_p\subset TM$ such that there exists a $1-$form $\alpha$, satisfying

\begin{enumerate}
\item $\ker_p(\alpha)=\xi_p$ for all $p \in M$ and
\item $\alpha\wedge d\alpha>0$.
\end{enumerate}
\end{defn}

\begin{figure}
     \centering
     \includegraphics{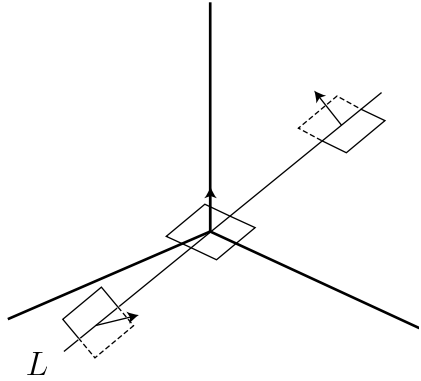}
     \caption{A rotationally symmetric contact structure on $\mathbb R^3$.}
     \label{basicexamplefig}
 \end{figure}

\begin{ex}\label{basicexample}
Figure~\ref{basicexamplefig} shows a family of planes in ${\mathbb R}^3$ that is invariant under rotation about the $z$-axis or translation in the $z$-direction.  The indicated line $L$ is {\rm Legendrian}, that is, $T_xL \subset \xi_x$ for all $x\in L$.  Note the planes twist slowly to the left as you move along L in either direction.  This example can be made more explicit by taking $\xi$ to be the kernel of $\alpha= rd\theta +dz$ and checking that it is a contact structure.
\end{ex}

The next definition and many of the results in the section are due to Giroux, \cite{Gi91}.

\begin{defn}
A surface $S\subset(M,\xi)$ is {\rm convex} if there exists  vector field $\vec{v}$ supported in a neighborhood of $S$ and transverse to $S$ such that flowing in the $\vec{v}$ direction preserves the contact planes.  If $\partial S\ne\emptyset$, we also require that $\partial S$ be Legendrian.  Such a $\vec{v}$ is called a {\rm contact vector field} for $\xi$.
\end{defn}

\begin{ex}\label{basic}
If $\xi$ is the contact structure of Example~\ref{basicexample}, it follows that any horizontal plane is  convex by considering the constant vector field $\vec{v}=\frac{\partial}{\partial z}$.  Indeed it follows that any surface in ${\mathbb R}^3$  transverse to the vector field $\frac{\partial}{\partial z}$ is convex.
\end{ex}

Roughly, $S$ is convex if and only if $S$ has a {\it product neighborhood}.  Convexity is a global condition; all smooth surfaces are locally convex.

\begin{defn} 
If $S\subset(M,\xi)$ is convex, {\rm the dividing set} is denoted $\Gamma_S$ and is defined to be $\{x\in S\mid \vec{v}(x) \in \xi_x \}$ .
\end{defn}

\smallskip
\noindent Intuition: If we think of the vector field $\vec{v}$ as vertical or perpendicular to $S$, then $\Gamma_S$ are those points whose contact planes are perpendicular to $S$.

\begin{figure}
    \centering
     \includegraphics{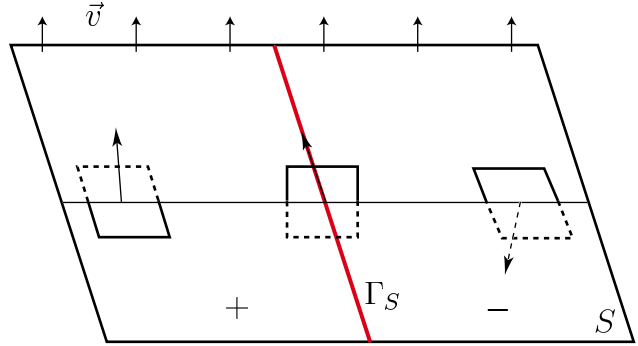}
     \caption{Dividing curve on a convex surface.}
     \label{localdividingcurve}
\end{figure}

\begin{defn}  The {\rm induced (singular) foliation} ${\mathcal F}_S$ on $S$ is defined by integrating the line field $\xi_p\cap T_pS$ on $S$.
\end{defn}

\begin{figure}
    \centering
     \includegraphics{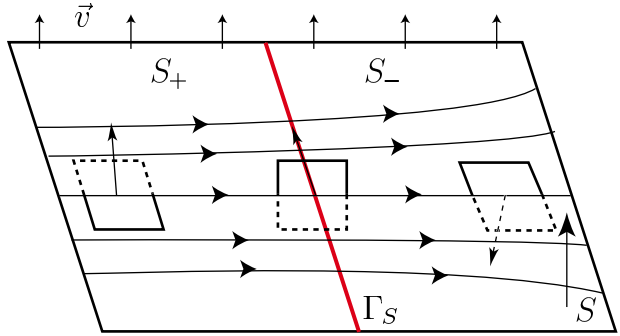}
     \caption{The induced foliation near a dividing curve.}
     \label{localinducedfol}
\end{figure}

Throughout this paper all manifolds and submanifolds are oriented.   Our contact structures are positive, thus the contact planes inherit a transverse orientation.   Orientations on $M$, $S$, and $\xi$  allow us to orient the leaves of ${\mathcal F}_S$.  Comparing these orientations leads to two ways of thinking about dividing sets:

\begin{enumerate}
\item The dividing set $\Gamma_S$ divides $S$ into regions where the contact planes are right side up or upside down relative to $S$ and as shown in Figure~\ref{localdividingcurve}.

\item With respect to the induced foliation on $S$, $\Gamma_S$ divides $S$ into {\it source} and {\it sink regions} as shown in Figure~\ref{localinducedfol}.
\end{enumerate}

\begin{figure}
    \centering
     \includegraphics{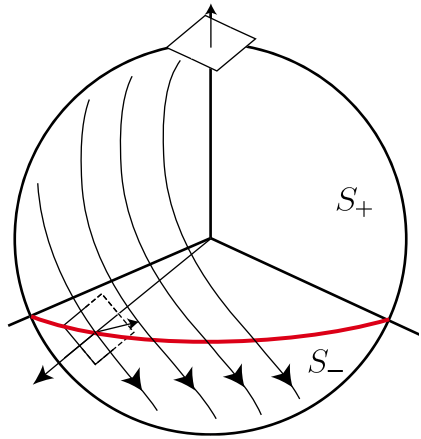}
     \caption{The induced foliation on a small sphere.}
     \label{sphere}
\end{figure}

\begin{ex}  Figure~\ref{sphere} shows the induced foliation and dividing set on a small round sphere about the origin of Example~\ref{basicexample}.  Contact structures are locally homogeneous by Pfaff's Theorem.  It follows that there exist small spheres like this about every point of any contact structure.
\end{ex}

\begin{prop}
The dividing set $\Gamma_S$ is a 1-dimensional submanifold of $S$ transverse to ${\mathcal F}_S$.
\end{prop}

\begin{proof} Choose coordinates $x\in S$ and $t$ in the $\vec{v}$ direction.
Then the 1-form defining $\xi$ may be written $\alpha=\beta(x)+f(x)dt$ where 
$\beta(x)$ is a 1-form on $S$, and $f$ is a function on $S$.  Since  $\ker\alpha=\xi$ we have:

\begin{enumerate}
\item $\alpha_x\left(\frac{\partial}{\partial t}\right) = 0$ if and only if  $f(x)=0$, and therefore
$\Gamma_S=f^{-1}(0)$.

\item $0\ne\alpha\wedge d\alpha =(\beta+fdt) \wedge (d\beta+df dt)
=\beta\wedge d\beta +\beta df dt +f dt d\beta = \beta df dt +f dt d\beta $. 
Therefore, if $f(x)=0$, then $\beta df dt \ne 0$, and in particular $df\ne 0$.  It now follows that $\Gamma_S=f^{-1}(0)$ is a submanifold of $S$.

\begin{figure}
    \centering
     \includegraphics{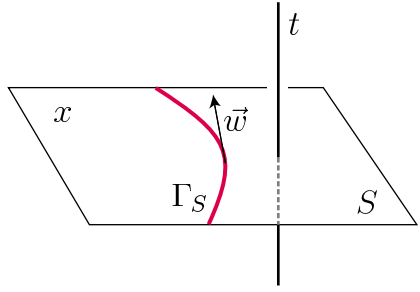}
     \caption{Coordinates near a convex surface.}
     \label{coords}
\end{figure}

\item Let $\vec{w}$ be tangent to $\Gamma_S$.  Then $dt(\vec{w})=df(\vec{w})=0$, so $\beta df dt \ne 0$ implies that $\beta(\vec{w}) \ne 0$, that is, $\vec{w} \notin ker\beta = T {\mathcal F}_S$.  It now follows that $\Gamma_S$ is transverse to the induced foliation on $S$.
\end{enumerate}
\end{proof}

The next several results will be used throughout the paper.

\begin{prop}  The isotopy class of $\Gamma_S$ does not depend on the choice of the contact vector field $\vec{v}$. \qed
\end{prop}

\begin{thm}[Existence of convex surfaces]\label{existence}
Every closed surface or compact surface with Legendrian boundary can be approximated by a convex surface.
\end{thm}

Theorem~\ref{existence} was proved by Giroux \cite{Gi91} for closed surfaces and by Honda for surfaces with boundary \cite{H1}.  Convex surfaces with Legendrian boundary were first used by Kanda~\cite{K97}.  Theorem~\ref{existence} follows from

\begin{prop}\label{sufficient} If $\partial S$ is Legendrian and the induced foliation, ${\mathcal F}_S$, is Morse-Smale, that is,

\begin{enumerate}
\item ${\mathcal F}_S$ has a finite number of closed leaves and Morse type singularities,
\item there are no saddle-saddle connections,
\item the holonomy about closed leaves is linear and either attracting or repelling,
\end{enumerate}
then $S$ is convex. \qed
\end{prop}

\begin{figure}
    \centering
     \includegraphics{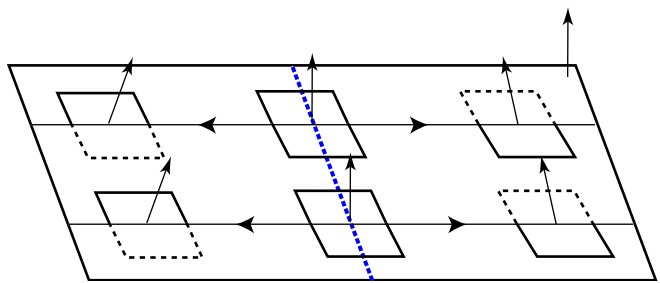}
     \caption{Legendrian divide.}
     \label{Legendriandivide}
\end{figure}

\begin{figure}
    \centering
     \includegraphics{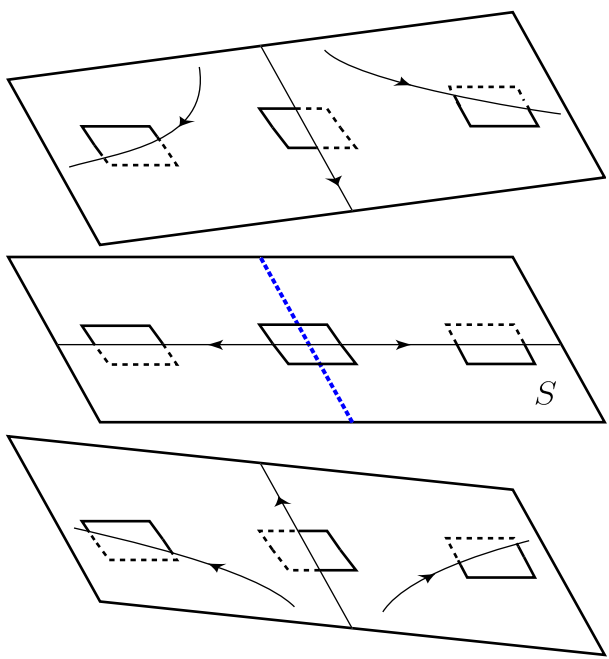}
     \caption{Flexibility near a Legendrian divide.}
     \label{flexibility}
\end{figure}

\begin{ex}\label{Legendrian} Proposition~\ref{sufficient} gives sufficient, but not necessary conditions for a surface to be convex.  Figure~\ref{Legendriandivide} shows a portion of convex surface whose induced foliation has a circle's worth of singularties.   The circle of singularities is called a {\rm Legendrian divide}.  The contact structure hinted at in the figure is invariant under translation in the vertical direction or parallel to the Legendrian divide.  Note that a Legendrian divide is not a dividing curve.

Figure~\ref{flexibility} shows a product neighborhood of the Legendrian divide and the effect of a slight perturbation of the original surface on the induced foliation.
\end{ex}

This example hints at a remarkable theorem about the possible induced foliations that can occur on perturbations of a convex surface.  Roughly, the Giroux Flexibility Theorem states that we can force the induced foliations to be whatever we like, within reason.  The statement of the theorem will make more sense after reading the definitions which follow it.

\begin{thm}[Giroux Flexibility Theorem \cite{Gi91}]
Let $S\subset(M,\xi)$ be a convex surface with dividing set $\Gamma_S$, and let 
${\mathcal F}$ be an arbitrary singular foliation on $S$ {\rm divided by} $\Gamma_S$, then there exists an isotopy of $S$ fixing $\Gamma_S$ (and keeping $S$ transverse to $\vec{v}$) such that at the end of the isotopy, $\mathcal F$ is the induced foliation on $S$.  
\end{thm}

\begin{defn}  $\Gamma_S$ {\rm divides} $\mathcal F$ if $\Gamma_S$ cuts $S$ into a maximal number of sink and source regions, that is, regions in which the induced foliation either points in at every boundary component or out at every component of each region.
\end{defn}

\begin{figure}
    \centering
     \includegraphics{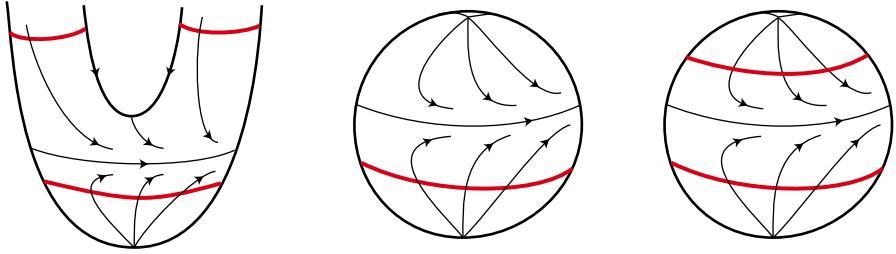}
     \caption{Hypothetical induced foliations.}
     \label{foliationsfig}
\end{figure}

\begin{ex}\label{foliations}
Figure~\ref{foliationsfig} shows three foliations.  The first and last are divided by the indicated curves, but the middle example is not; it is not cut into a maximal number of sink and source regions as the third example shows.  Another example  to which the Giroux Flexibility Theorem can be applied is to replace the indicated foliation on the annular region between the two dividing curves of the third example with a Legendrian divide.
\end{ex}

\begin{defn}  $(M,\xi)$ is {\rm tight} if there does not exist an embedded disk $D\subset M$ such that $D$ is tangent to $\xi$ along its boundary (i.\,e.,\  $T_xD=\xi_x$ for all $x\in\partial D$).  $(M, \xi)$ is called {\rm overtwisted} if it is not tight.
\end{defn}

Overtwisted contact structures are classified by their underlying 2-plane bundles \cite{E89}.  The notion of tightness is analogous to tautness or non-existence of Reeb components in foliation theory or incompressibility of surfaces.  We shall see that tight contact structures reflect the underlying topology of the 3-manifolds which carry them.  

\begin{figure}
    \centering
     \includegraphics{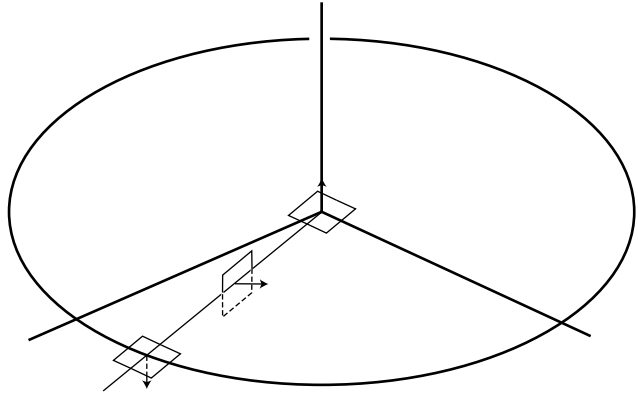}
     \caption{An overtwisted disk.}
     \label{overtwisted}
\end{figure}

Figure~\ref{overtwisted} shows an overtwisted disk that would live in the contact structure described in Example~\ref{basicexample} if the contact planes were allowed to rotate too quickly along rays leaving the origin.

\begin{prop}[Giroux \cite{Gi01}]\label{tightneighborhood}  If $S\subset(M,\xi)$ is convex, a product neighborhood of $S$ is tight if and only if one of the following is satisfied:
 \begin{enumerate}
\item $S=S^2$, and $\Gamma_S$ is connected

\item $S\ne S^2$, and {\it no component of $\Gamma_S$ is null-homotopic in $S$}.
\end{enumerate}
\end{prop}

\begin{proof}[Sketch]
$(\Rightarrow)$  If either (1) or (2) is false, use the Giroux Flexibility Theorem to realize a null-homotopic Legendrian divide, as discussed in Example~\ref{foliations}.  The disk in $S$ bounded by the Legendrian divide will be an overtwisted disk.

$(\Leftarrow)$ We need a starting point and gluing theorems.  That is, until this point, we have not even stated that there are any tight contact structures on any manifold.  The next theorem will address this.  Given simple examples of tight contact structures we require gluing theorems to produce more complicated examples.  This paper will eventually describe several gluing theorems.  Another strategy, used by Giroux, is to produce models in which the desired $S$ and $\Gamma_S$ exist and must be tight.
\end{proof}

\begin{thm}\label{B3}
There exist a tight contact structure on $B^3$, moreover, two tight contact structures which induced the same foliations on $B^3$ are diffeomorphic.  \qed
\end{thm}

The existence portion of the theorem is due to Bennequin \cite{Be83}, and the uniqueness portion is due to Eliashberg \cite{E92}.  In light of the Giroux Flexibility Theorem, we will paraphrase Theorem~\ref{B3} by saying that there is a unique tight contact structure on $B^3$.

Convex surfaces are required to have Legendrian boundary.  Therefore to decompose manifolds with convex boundaries along convex surfaces, we will need to know which curves on a convex surface $S$ can be ``made Legendrian''.  That is, we need to know which curves are contained in the leaves of some foliation $\mathcal F$ divided by $\Gamma_S$.   $S$ can be perturbed so that they become Legendrian.   The next definition and theorem of Honda's \cite{H1} exactly answers this question.

\begin{defn} 
A properly embedded 1-submanifold $C$  of a convex surface $S$ is {\rm non-isolating} if
\begin{enumerate}

\item $C$ is transverse to $\Gamma_S$ and

\item the closure of every component of $S\backslash (\Gamma_S\cup C)$ intersects $\Gamma_S$.
\end{enumerate}
\end{defn}

\begin{thm}[Legendrian Realization Principle]
If $C$ is non-isolating then $C$ can be made Legendrian.
\end{thm}

\begin{proof}[Sketch]
The non-isolating condition guarantees that $C$ can be extended to a foliation divided by $\Gamma_S$.  Then use Giroux Flexibility to realize this foliation on $S$.
\end{proof}

\begin{figure}
    \centering
     \includegraphics{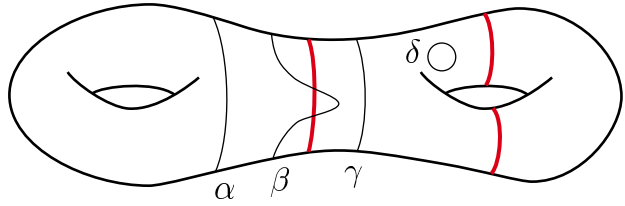}
     \caption{Hypothetical Legendrian curves.}
     \label{hypo}
\end{figure}

\begin{ex}  Of the curves shown in Figure~\ref{hypo}, only $\beta$ and $\gamma$ are non-isolating.  Notice that any curve, such as $\beta$, which intersects $\Gamma_S$ is non-isolating.  It is not too hard to extend, say $\beta$, to a singular foliation on $S$ divided by $\Gamma_S$, however $\beta$ will end up passing through singularities, that is, it will not be a smooth curve on $S$.  Note also that in the definition of non-isolating, $C$ is not necessarily connected or closed.
\end{ex}

\section{Preview}

At this point we have enough of the foundational tools in place to sketch, in general terms, some of the issues and techniques involved in studying contact structures from a cut-and-paste point of view.

\medskip
\noindent {\bf Classification:}  Given a 3-manifold $M$ and a collection of curves $\Gamma$ contained in $\partial M$, how many tight contact structures, up to equivalence, are there on $M$ with $\Gamma_{\partial M}=\Gamma$?  Equivalence might be either diffeomorphism or isotopy taking one contact structure to another.  
\smallskip

\begin{figure}
    \centering
     \includegraphics{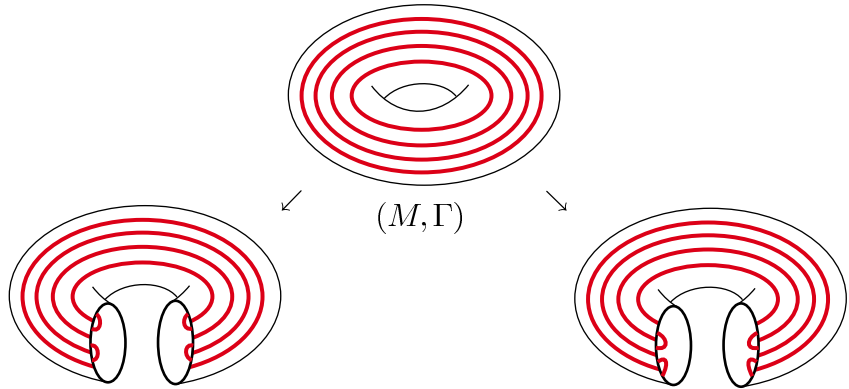}
     \caption{Different convex decompositions.}
     \label{preview}
\end{figure}

To be specific, consider the case of a solid torus with four dividing curves on its boundary shown in Figure~\ref{preview}.

\medskip
\noindent {\bf Decomposition:}  How many ``sensible'' ways are there to decompose such an $(M, \Gamma)$?
\smallskip

Continuing with the solid torus example, Figure~\ref{preview} suggests that there are just two possible decompositions, thus, there are at most two tight contact structures carried by $(M, \Gamma)$.

\medskip
\noindent {\bf Gluing:} Which of the decompositions into tight pieces can be glued to form a tight union?
\smallskip

Unlike many situations in 3-dimensional topology, it is very difficult to give general conditions under which the union of tight pieces is tight.  The problem is that a manifold can contain a large overtwisted disk, but when it is chopped into small pieces, none of the pieces may contain overtwisted disks themselves.

It turns out that regluing either of the decompositions shown in Figure~\ref{preview} gives a tight contact structure.  {\it A priori}, we do not know that these two contact structures are different.  By gluing we can conclude only that $(M, \Gamma)$ carries at least one tight contact structure.

\medskip
\noindent {\bf Invariants:} Of the various ways of gluing into a tight union, which result in non-isotopic contact structures?
\smallskip

In our example, an Euler charactistic type invariant shows that the two gluings result in different contact structures.  It follows that $(M, \Gamma)$ carries exactly two tight contact structures.

\section{Convex Decompositions}

A convex decomposition can be viewed in two ways.  First, you can start with a contact structure on a 3-manifold $M$ and keep splitting $M$ along convex surfaces until the pieces are balls.   Alternatively, you can start with $M$ and a collection of curves $\Gamma$ on $\partial M$ that you hope will end up being dividing curves for a contact structure that you are trying to build, and then split along surfaces which you hope will end up being convex.  We need to see how actual convex surfaces intersect so that this structure can be correctly modelled in the definition of a convex decomposition.

\begin{figure}
    \centering
     \includegraphics{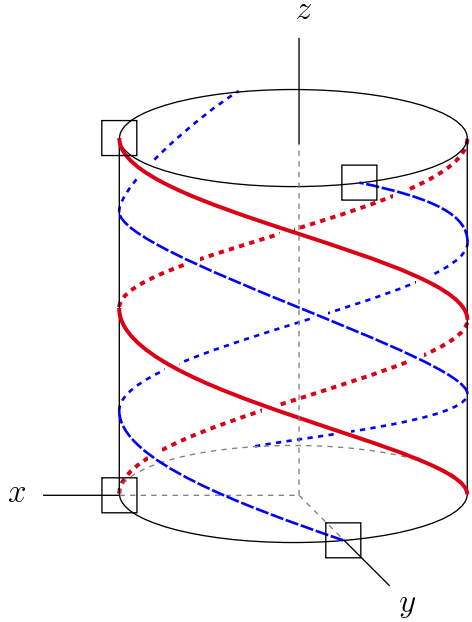}
     \caption{The neighborhood of a Legendrian curve.}
     \label{Legendriancurves}
\end{figure}

\begin{ex}\label{Legendrianex} The kernel of $\alpha_k=\sin(2\pi kz) dx +\cos(2\pi kz)dy$ defines a contact structure on $\mathbb R^3$ shown in Figure~\ref{Legendriancurves}.  In this example, the contact planes all contain the $z$-axis, that is, any vertical line is Legendrian.  The foliation induced on horizontal planes is a linear foliation with slope changing as the height of the plane increases.  The vector field given by $\frac{\partial}{\partial r}$ in cylindrical coordinates is a contact vector field, thus a cylinder at constant distance from the $z$-axis is convex.  The dividing curves on this cylinder start on the $x$-axis and spiral upwards at a rate depending on $k$.  Figure~\ref{Legendriancurves} also shows the tangencies of the contact planes and the cylinder as long dashed lines starting on the $y$-axis,.

By restricting to the cylinder about the $z$-axis and identifying top and bottom, $\alpha_k$ also defines a contact structure on $S^1 \times D^2$.  The key features of this contact structure are:
\begin{enumerate}

\item $T=\partial(S^1 \times D^2)$ is convex

\item $\#\Gamma_{T}=2$

\item $\mbox{slope}(\Gamma_{T})=-\frac1k$
\end{enumerate}

\end{ex}

The next theorem may be paraphrased by saying that Legendrian curves, such as the quotient of the $z$-axis in the previous example, have standard neighborhoods.
\begin{thm}[Kanda \cite{K97}, Makar-Limanov \cite{ML98}]
There is a unique tight contact structure on $S^1 \times D^2$ such that (1), (2), and (3) hold. \qed
\end{thm}

When we start with a manifold with convex boundary and cut it along a convex surface, the cutting surface, by definition of convexity, intersects the boundary in a Legendrian curve.  The next example is a portion of the region shown in Figure~\ref{Legendriancurves}.  From it we see how the dividing curves on a pair of intersecting surfaces are related near their Legendrian curve of intersection.

\begin{figure}
    \centering
     \includegraphics{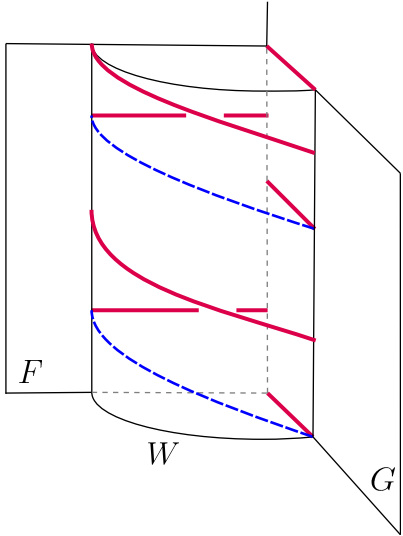}
     \caption{Intersecting convex surfaces.}
     \label{convexintersectionfig}
\end{figure}

\begin{ex}\label{convexintersection}
In Example~\ref{Legendrianex} the $xz$-plane is convex with respect to the contact vector field
$\frac{\partial}{\partial y}$ and similarly the $yz$-plane is convex with respect to $\frac{\partial}{\partial x}$.  Figure~\ref{convexintersectionfig} shows portions of these planes, labelled $F$ and $G$ and their dividing curves.   Notice that $\Gamma_F$ and $\Gamma_G$ are horizontal lines starting on the $z$-axis and ending at a point of tangency of a contact plane and the vertical cylinder.
\end{ex}

From this we see that for general intersecting convex surfaces $F$ and $G$, the endpoints of $\Gamma_F$ and $\Gamma_G$ alternate along curves of $F\cap G$.  Further examination of Figure~\ref{convexintersectionfig} shows that if the corner of the wedge $W$ subtended by $F$ and $G$ is smoothed, the manifold produced has convex boundary and the dividing curves of $F$ and $G$ are joined by turning to the right (when viewed from the outside of $W$).  The ``turn to the right'' rule that is forced on us in the presence of a {\it positive} contact structure serves as the model for defining the orientation conventions in convex decompositions.

\begin{figure}
    \centering
     \includegraphics{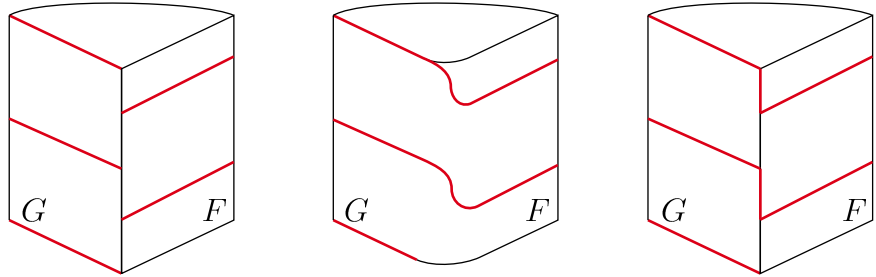}
     \caption{Dividing curves before and after smoothing.}
     \label{round}
\end{figure}

Figure~\ref{round} shows three views of $W$.  The first shows $W$ before rounding corners, the second is after rounding corners.  The last picture shows $W$ without the corner rounded, but it shows the effect on the dividing curves of corner rounding.  Most of the figures in this paper are drawn in this fashion.

\begin{figure}
    \centering
     \includegraphics{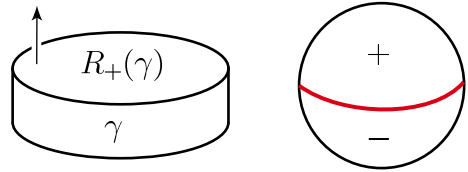}
     \caption{$B^3$ with a single suture.}
     \label{suture}
\end{figure}

\begin{defn}
$(M,\gamma)$ is a {\rm sutured manifold} if:
\begin{enumerate}

\item $\gamma \subset\partial M$ is a union annuli and tori,

\item $(\partial M) \backslash\gamma$ is a disjoint union of two subsurfaces $R_+(\gamma)$ and $R_-(\gamma)$, and

\item crossing an annular suture takes you from $R_\pm(\gamma)$ to $R_\mp(\gamma)$.

\end{enumerate}
\end{defn}
Gabai \cite{Ga} defined sutured manifolds to study taut foliations.  We will primarily be concerned with that case that all sutures are annuli.  Figure~\ref{suture} shows two views of a solid ball with a single annular suture.  The first view shows a manifold with corners, as it should be drawn.  The second shows how sutures will be drawn; the manifold appears smooth, and the sutures are very skinny.

\begin{figure}
    \centering
     \includegraphics{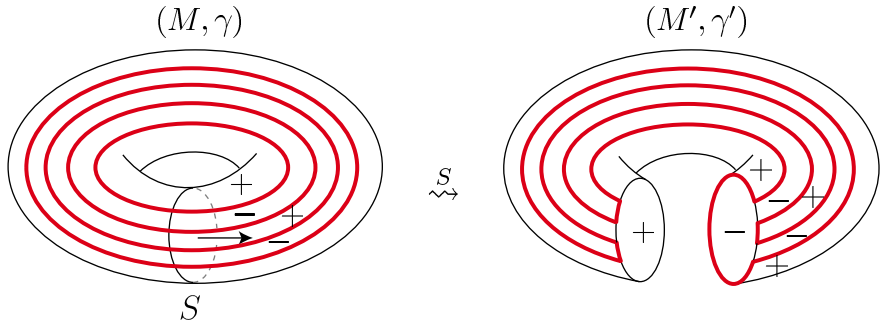}
     \caption{Sutured manifold splitting.}
     \label{sutureddecomposition}
\end{figure}

\begin{defn}
If $S$ is an oriented properly embedded surface in $(M, \gamma)$, $(M,\gamma)\stackrel{S}{\rightsquigarrow}(M',\gamma')$ is defined by $M'=M\backslash S$ and introducing sutures as needed to separate the positively and negatively oriented portions of $\partial(M\backslash S)$ as shown in Figure~\ref{sutureddecomposition}.
\end{defn}

\begin{figure}
    \centering
     \includegraphics{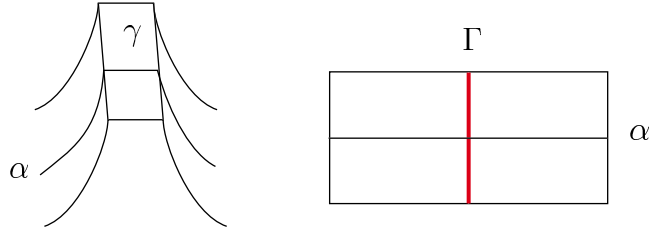}
     \caption{An arc crossing a suture compared to an arc crossing a dividing curve.}
     \label{suturevsconvex}
\end{figure}

\begin{defn}  A {\rm convex structure} is a pair $(M,\Gamma)$ such that:
\begin{enumerate}

\item $\Gamma$ is a disjoint union of curves in $\partial M$,

\item $\partial M$ split along $\Gamma$ is the disjoint union of two subsurfaces, $R_+(\Gamma)$ and $R_-(\Gamma)$, and

\item crossing a dividing curve takes you from $R_\pm(\Gamma)$ to $R_\mp(\Gamma)$.

\end{enumerate}

Also assume each component of $\partial M$ has dividing curves on it.
\end{defn}
A sutured manifold $(M,\gamma)$ is a manifold with corners.  Convex structures $(M,\gamma)$ are smooth.  Figure~\ref{suturevsconvex} shows portions of a sutured manifold near a suture and a convex structure near a dividing curve.   Notice that the 2-planes along each arc $\alpha$ turn over as $\alpha$ is traversed, but they do so in different fashions.

Of course we hope that $(M, \Gamma)$ will carry an actual tight contact structure, just as Gabai would like $(M,\gamma)$ to carry a taut foliation; however there is no {\it a priori} reason that it will.  When discussing a surface with curves on it, such as $(S, \sigma)$ in the next definition, there is no need to distinguish between an ``abstract'' convex surface and an actual convex surface, for a contact structure is uniquely determined in a (product) neighborhood of $S$ by the dividing curve configuration $\sigma$.

\begin{figure}
    \centering
     \includegraphics{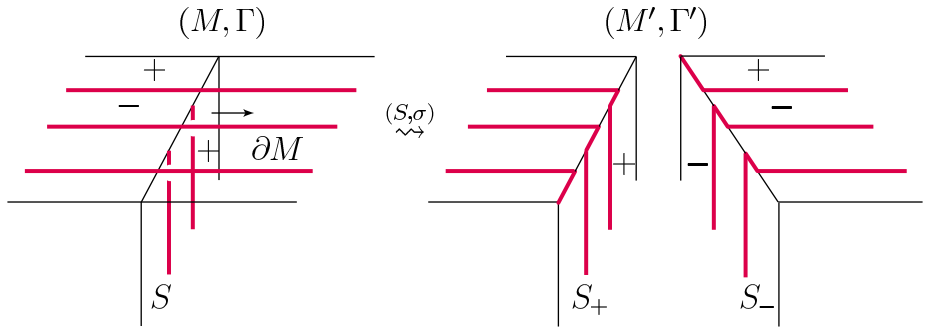}
     \caption{A convex splitting,}
     \label{convexdecomp}
\end{figure}

\begin{defn} Let $(S,\sigma)$ be a convex surface in $(M,\Gamma)$ such that $\partial S$ is non-isolating in  $\partial M$, and the endpoints of $\sigma$ alternate with points of $\Gamma \cap \partial S$ along $\partial S$.  Define $(M,\Gamma)\stackrel{(S,\gamma)}{\rightsquigarrow}(M',\Gamma')$ by $M'=M\backslash S$ and by adding new portions of dividing curves to $(\sigma \cup \Gamma)\backslash S$ using the ``turn to the right'' rule shown in Figure~\ref{convexdecomp}.
\end{defn}

A sutured manifold with annular sutures $(M, \gamma)$ naturally determines a convex structure $(M, \Gamma)$ by replacing each annulus of $\gamma$ by its core.  To be able to use Gabai's existence theorems for sutured manifold decompositions in our setting we must be able to start with a sutured manifold splitting (the top row of the following diagram) and then produce a convex surface, $(S, \sigma)$, such that the diagram commutes.  We discuss how this can be done through a series of examples.

\[  \begin{array}{ccc}
(M,\gamma)&\stackrel{S}{\rightsquigarrow}&(M',\gamma')\\
\downarrow &  & \downarrow\\
(M,\Gamma)&\stackrel{(S,\sigma)}{\rightsquigarrow}&(M',\Gamma')
\end{array} \]

\begin{figure}
    \centering
     \includegraphics{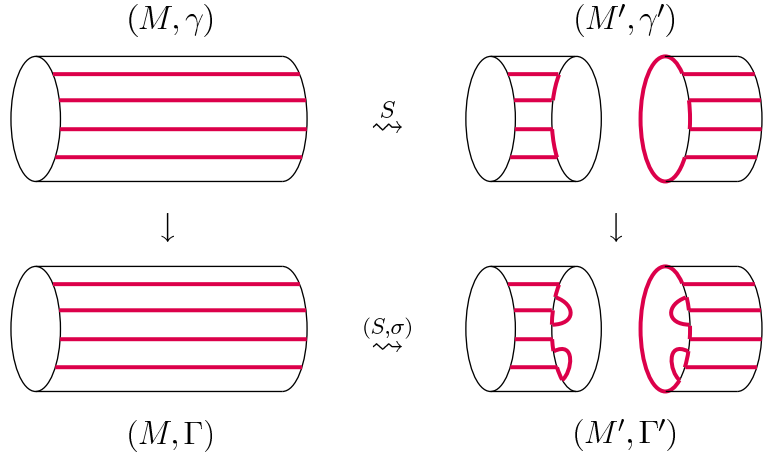}
     \caption{A commuting diagram of sutured manifold and convex splittings.}
     \label{commutative}
\end{figure}

\begin{ex}  Figure~\ref{commutative} shows how to introduce {\rm boundary-parallel} dividing curves, $\sigma$, so that the diagram commutes.  This technique works, provided that every component of $\partial S$ has nonempty intersection with $\Gamma$.
\end{ex}
 
\begin{defn}\label{boundaryparallel}
A convex surface $(S, \sigma)$ has {\rm boundary-parallel dividing curves} if $\partial S$ is nonempty, every component of $\partial S$ intersects $\sigma$, and $\sigma$ is collection of arcs each of which bounds a half disk that contains a portion of $\partial S$ but no other arcs of $\sigma$.
\end{defn}

\begin{figure}
    \centering
     \includegraphics{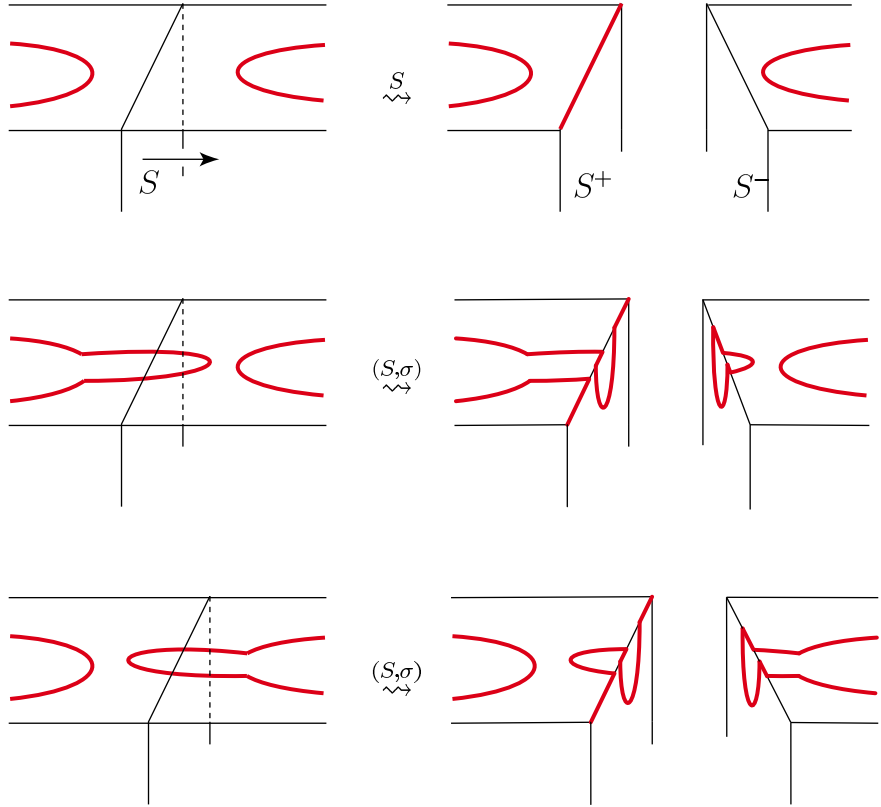}
     \caption{Two possibilities for $(S, \sigma)$.}
     \label{sutdecomp}
\end{figure}

\begin{ex}\label{sutdecomposition} Now consider the possibility that $\partial S \cap \Gamma = \emptyset$.  Such an $S$ might have isolating boundary, that is, it might not be possible to make it Legendrian and hence $S$ convex.  Figure~\ref{sutdecomp} shows first a sutured manifold splitting along a surface $S$ with $\partial S \cap \Gamma = \emptyset$.  The second two portions of the figure show two possible ways of introducing intersections between $\partial S$ and $\Gamma$ and of adding boundary compressible $\sigma$ to $S$.  

There are two key features in this example.  First, the strategy of introducing intersections can only work if there are dividing curves on the same component of $\partial M$ as $\partial S$ -- this will show up in the definition of ``sutured manifold with annular sutures'' below.  And second, only one of the perturbations of $S$ makes the splitting diagram commute.  
\end{ex}

\begin{figure}
    \centering
     \includegraphics{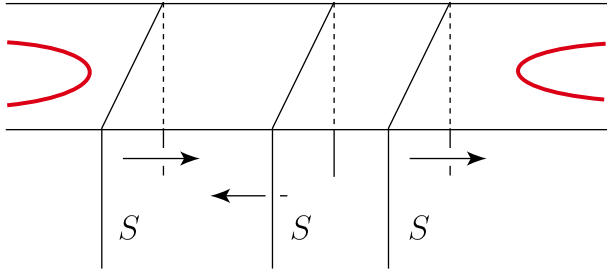}
     \caption{There is no consistent way to introduce intersections between $\partial S$ and $\Gamma$.}
     \label{wellgroomed}
\end{figure}

Figure~\ref{wellgroomed} is similar to Figure~\ref{sutdecomp} in that $\partial S \cap \Gamma = \emptyset$, but in this case there are multiple portions of $\partial S$, each with its own orientation preference for creating a pair of intersections with $\Gamma$,  and they can not all be satisfied simultaneously.  Rather than describe how to get around this, we just point out that Gabai confronted a similar situation in developing sutured manifold theory.  He introduced a notion of ``well-groomed'' sutured decompositions, that is, he showed that splittings could be assumed to have coherently oriented boundary components, and for such splittings we can produce a commutative splitting diagram using the technique of Example~\ref{sutdecomposition}.

\begin{thm}\label{equivalent} Let $(M,\gamma)$ be an irreducible {\rm sutured manifold with annular sutures}, and let $(M,\Gamma)$ be the corresponding convex structure.  The following are equivalent.

\begin{enumerate}
\item $(M,\gamma)$ is {\rm taut}.
\item $(M,\gamma)$ has a {\rm sutured manifold decomposition}.
\item $(M,\gamma)$ carries a {\rm taut foliation}.
\item $(M,\gamma)$ carries a {\rm universally tight} contact structure.
\item $(M,\gamma)$ carries a tight contact structure.
\end{enumerate}
\end{thm}

\begin{defn}
A sutured manifold has {\rm annular sutures} if $\partial M$ is nonempty, every boundary component contains at least one annular suture, and if there are no toroidal sutures.

$(M,\gamma)$ is {\rm taut} if $R_+(\gamma)$ and $R_-(\gamma)$ are Thurston norm minimizing in their homology class in $H_2(M,\gamma)$.

A {\rm sutured manifold decomposition} of $M$ is a sequence of splittings
$$
(M,\gamma)\stackrel{S_1}{\rightsquigarrow}
\cdots\stackrel{S_m}{\rightsquigarrow}
\cup(B^3,S^1 \times I)
$$
where $(B^3,S^1 \times I)$ denotes the sutured manifold shown in Figure~\ref{suture}.

A foliation is {\rm taut} if every leaf intersects a closed transversal.

A contact structure is {\rm universally tight} if $(\widetilde M,\tilde{\xi})$ is tight.
\end{defn}

Thurston \cite{Th} proved (3) implies (1).  Gabai \cite{Ga} proved (1) implies (2) and (2) implies (3).  Eliashberg and Thurston \cite{ET} showed (3) implies (4).  It is immediate that (4) implies (5).  All of these results apply without the additional assumption of annular sutures.  Since $S^3$ carries a tight contact structure but can not support a taut foliation, some additional hypothesis is necessary for (5) to imply (1).  

The techniques of (5) implies (1) will not be used in the rest of the paper, so we will instead sketch a direct proof of (2) implies (4) that has the advantages of making the importance and utility of universal tightness clear. The proof introduces a gluing strategy that will be used repeatedly.

\begin{proof}[Proof of (2) implies (4)]
First replace the given sutured manifold decomposition with a corresponding convex decomposition
$$
(M,\Gamma)\stackrel{(S_1,\sigma_1)}{\rightsquigarrow} \cdots\stackrel{(S_m, \sigma_m)}{\rightsquigarrow}
\cup(B^3,S^1)
$$

By Theorem~\ref{B3}, $(B^3, S^1)$ carries a (universally) tight contact structure.  By construction, the surfaces $(S_i, \sigma_i)$ have boundary-parallel dividing curves (see Definition~\ref{boundaryparallel}) thus this portion of the theorem follows from the next gluing theorem.
\end{proof}

\begin{thm}\label{Colin}(Colin \cite{Colin2})
Let $(M,\Gamma)\stackrel{(S,\sigma)}{\rightsquigarrow}(M',\Gamma')$.  If $M$ is irreducible, $S$ has boundary compressible dividing curves, and $(M',\Gamma')$ carries a universally tight contact structure then so does $(M,\Gamma)$.
\end{thm}

\begin{proof}[Sketch] We will illustrate key ideas of our interpretation \cite{HKM2} of Colin's gluing theorem with examples.  The proof strategy is to:
\begin{enumerate}
\item Suppose the contact structure on $(M,\Gamma)$ obtained by gluing $(M', \Gamma')$ along $(S, \sigma)$ is overtwisted, and let $D$ be an overtwisted disk.

\item In {\it small steps} isotop $S$ to $S'$ and eventually off $D$.

\item While isotoping $S$, make sure that $M\backslash S'$ stays universally tight. 

\end{enumerate}

This strategy gives a contradiction once $S'\cap D=\emptyset$, for $M\backslash S'$ is both tight and contains the overtwisted disk $D$.

In {\it small steps} refers to a fundamental  idea due to Honda \cite{H3}.  That is, any isotopy of a convex surface $S$ can be expressed as a sequence of {\it bypasses} or their inverses.

\begin{defn}
A {\rm bypass} consists of:
\begin{enumerate}

\item a Legendrian arc $\alpha$ connecting 3 dividing curves in $S$.

\item a Legendrian arc $\beta$ joining $\partial\alpha$

\item a convex half disk in $M\backslash S$ with boundary equal to $\alpha \cup \beta$ which contains a single dividing curve.

\end{enumerate}
\end{defn}

\begin{figure}
    \centering
     \includegraphics{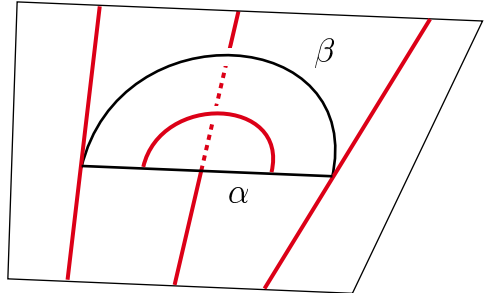}
     \caption{A bypass attached to $S$ along $\alpha$.}
     \label{bypass}
\end{figure}

Figure~\ref{bypassmove} shows the effect on $\sigma$ of isotoping $S$ across a bypass.  Notice that the ``turn to the right'' rule for dividing curves going around corners looks more like a ``turn to the left'' rule when it is viewed from inside the manifold.

\begin{figure}
    \centering
     \includegraphics{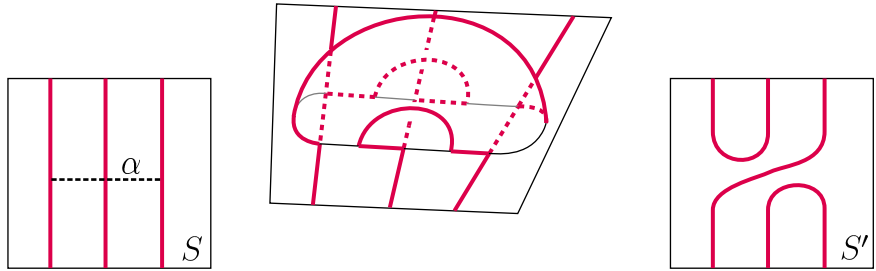}
     \caption{The effect on $\Gamma_S$ of isotoping $S$ through a bypass attached along $\alpha$.}
     \label{bypassmove}
\end{figure}

\begin{ex}\label{parallel}
The global effect on the dividing curves of a bypass move depends very much on how the local picture sits with respect to the entire surface and dividing curve set.  Figure~\ref{parallelfig} shows an example in which the arc of attachment connects two parallel dividing curves to a third.  Isotoping $S$ across this bypass has the effect of removing two dividing curves from $S$.
\end{ex}

\begin{figure}
    \centering
     \includegraphics{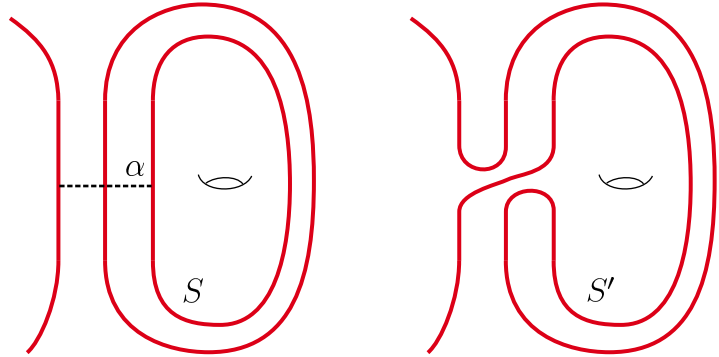}
     \caption{This bypass removes two parallel dividing curves.}
     \label{parallelfig}
\end{figure}

Continuing with the gluing theorem, we now consider some examples of bypasses that $S$ might have to be isotoped through while moving $S$ off of $D$.  Hopefully universal tightness of $M\backslash S'$ follows from universal tightness of $M\backslash S$ in each case.

\begin{ex}
Notice that the dividing curves in Figure~\ref{gluing1} are boundary-parallel.  If such a bypass were to exist in $M$, then, as shown, $S'$ would contain a null-homotopic dividing curve.  Proposition~\ref{tightneighborhood} implies the existence of an overtwisted disk near $S'$.  Since $S'$ may be thought of as living in the complement of $S$, and we are assuming $M\backslash S$ is universally tight, the bypass drawn in Figure~\ref{gluing1} can not exist.
\end{ex}

\begin{figure}
    \centering
     \includegraphics{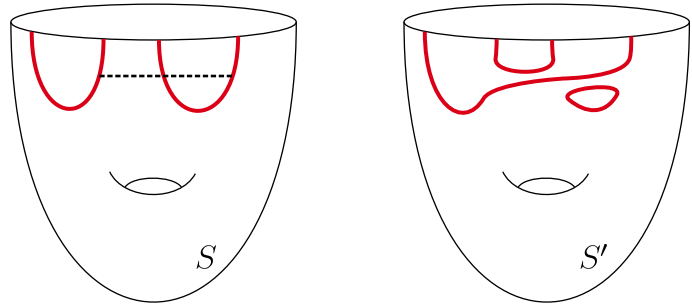}
     \caption{This bypass can not exist in a tight contact structure.}
     \label{gluing1}
\end{figure}

\begin{ex}\label{trivial}
Figure~\ref{gluing2} shows a similar-looking bypass that has a very different effect on the dividing curves of $\sigma$.  Up to isotopy, the dividing curves are unchanged.  Though we omit the proof here (see \cite{HKM2}), it is a consequence of the uniqueness of tight contact structures on a ball that $S$ and $S'$ cobound a contact product.  It follows that $M\backslash S$ and $M\backslash S'$ are contactomorphic, hence both are universally tight.
\end{ex}

\begin{figure}
    \centering
     \includegraphics{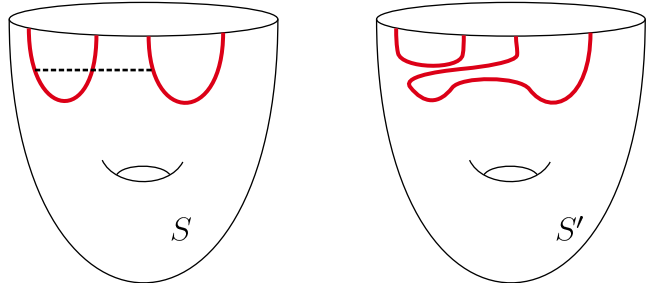}
     \caption{A trivial bypass.}
     \label{gluing2}
\end{figure}

\begin{ex}
Figure~\ref{gluing3} is the most mysterious of these examples.  Such a bypass may or may not exist, and there is no reason for $M\backslash S'$ to be tight if we only assume tightness of $M\backslash S$.  However,  $\alpha$ represents a nontrivial element of $\pi_1(S)$ or $\pi_1(M)$, and any cover $\widetilde{M\backslash S}$ is still tight (by universal tightness).  If we lift to the right cover, $\alpha$ is unwound, and this example becomes the same as the previous example.   Thus the proof strategy may be continued in this and subsequent covers.
\end{ex}
\begin{figure}
    \centering
     \includegraphics{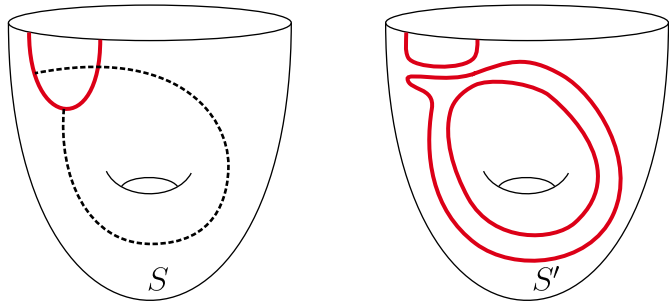}
     \caption{It isn't clear what the implications of such a bypass are in general.}
     \label{gluing3}
\end{figure}

There are several other types of bypass configurations to check, but this pattern repeats itself.  Bypasses are of three types.  Those which can not exist.  Those which cause no trouble if they do exist, and all of the rest.  The typical situation is that troublesome bypasses can be dealt with in the right cover.  This is the point and power of the assumption of universal tightness.
\end{proof}

\section{Tori}

\begin{figure}
    \centering
     \includegraphics{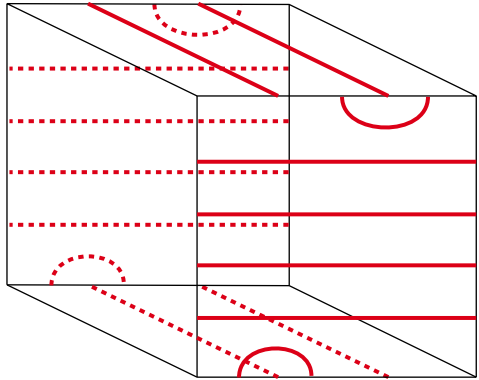}
     \caption{$\xi_2$ on $T \times I$.}
     \label{torustorsion}
\end{figure}

\begin{ex}\label{torus}  Let $\alpha_k=\sin(2k\pi z)dx +\cos(2k\pi z)dy$ as in Example~\ref{Legendrianex}.  Restricting $\alpha_k$ to the cube $[0,1] \times [0,1] \times [0,1]$ and identifying the front with the back face and the left with the right face defines a contact structure on $T \times I$.  Neither $T \times \{0\}$ nor $T \times \{1\}$ is convex.  Perturbing $T \times \{0\}$ and $T \times \{1\}$ so that they are convex gives the contact structure $\xi_k$ on $T \times I$ shown in Figure~\ref{torustorsion}.    A vertical annulus, such as the one shown on the front face, will be convex and  have $2k$ closed dividing curves.
\end{ex}

The next theorem gives a very general, but rough, classification theorem for tight contact structures.

\begin{thm}
If $M$ is irreducible, then $M$ carries finitely many tight (or universally tight) contact structures if and only if $M$ is atoroidal.
\end{thm}

The if direction is due to Colin, Honda, and Giroux \cite{CHG}, and the only if direction is due to Colin \cite{Co99b, Co01} and Honda, Kazez, and Mati\'c \cite{HKM2}.

\begin{proof}
We will explain the following portions of the proof of finiteness direction:

\begin{enumerate}

\item There is a finite collection of branched surfaces in $M$ which carry every tight contact structure.

\item If a branched surface carries infinitely many tight contact structures then it carries tori.
\end{enumerate}

\smallskip
\noindent Proof of (1).

\begin{itemize}

\item Pick a triangulation $\tau$ of $M$, and isotop it until $\tau^1$ is a collection of Legendrian arcs.

\item Isotop $\tau^2$ relative to $\tau^1$ so that each face is convex.

\item Isotop $\tau$ to remove interior $\partial$-parallel dividing curves.  Figure~\ref{simplicialbypass} shows how to pry a two cell open along an edge to effect a bypass move and accomplish this.

\begin{figure}
    \centering
     \includegraphics{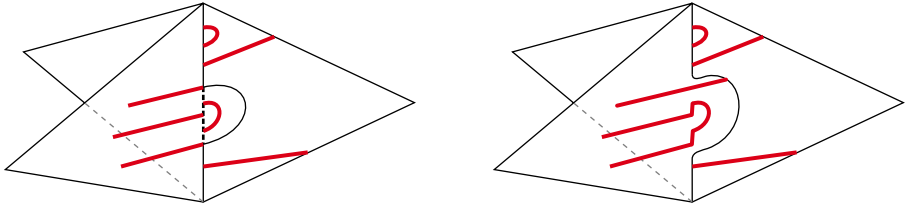}
     \caption{Isotoping an edge of $\tau$ across a bypass.}
     \label{simplicialbypass}
\end{figure}

\item For each $\Delta \in \tau^3$, group the dividing curves in $\partial \Delta^3$ so that, except for a bounded number of dividing curves near the vertices,  they are contained in at most $5$ prisms $P_i$.  See Figure~\ref{tetra}.

\begin{figure}
    \centering
     \includegraphics{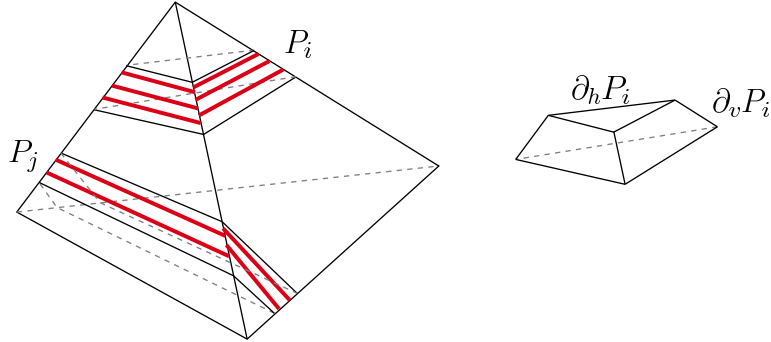}
     \caption{Dividing curves on $\partial \Delta$ are carried by a family of prisms.}
     \label{tetra}
\end{figure}

\item Use Giroux Flexibility to force the foliation induced by $\xi$ on each $\partial_v P_i$ to be a union of vertical arcs and on $\partial_hP_i$ to be a fixed non-singular foliation.

\item By the uniqueness of tight contact structures on $B^3$ we may assume all vertical arcs in $P_i$ are Legendrian.

\item The union over $\Delta \in \tau^3$ is naturally a neighborhood $N(B)$ of a branched surface $B$, and by construction there are only finitely many such $B$.

\item Each component of $\Delta^3-N(B)$ is a polygonal ball.  The number of such polygonal balls is bounded, and the possible dividing curve configurations on the boundary faces of each polygonal ball is also bounded.  Thus $\xi$ is defined by $\xi\mid N(B)$ up to finitely many choices.
\end{itemize}

\smallskip
\noindent Proof of (2).

Suppose two contact structures $\xi_0,\xi_1$ are carried by $B$.  By construction the foliations induced on $\partial_h N(B)$ agree, thus $\xi_1$ is defined by $\xi_0$ and a finite set of integer weights on the sectors of $B$ which describe the twisting of the planes of $\xi_1$ relative to the planes of $\xi_0$ along vertical Legendrian arcs of $N(B)$.

An infinite collection of contact structures all carried by one branched surface give an infinite collection of weights.   Since the contact structures are all positive, there is a lower bound, perhaps negative, on these weights.  It follows that there must be a non-negative collection of integer weights on $B$.  In the standard way, these non-negative weights can be used to piece together a surface in $N(B)$ that is transverse to the vertical Legendrian arcs.  The induced foliation on such a surface has no singularities, thus the surface is either a torus or a Klein bottle.

We draw two conclusions from this portion of the argument.

\begin{itemize}
\item Changing weights along a torus doesn't change the homotopy class of the 2-plane bundle, thus it follows that only finitely many 2-plane bundles support tight contact structures.

\item The only way to produce infinitely many contact structures on a given space is to insert $\xi_k$, as defined in Example~\ref{torus}, in a neighborhood of a torus.
\end{itemize}

With this in mind we sketch some of the remaining steps in the infinitely many portion of the theorem.

\begin{enumerate}
\setcounter{enumi}{2} 
\item A toroidal manifold has a universally tight contact structure.
\item Inserting $\xi_{k}$ near the torus preserves universal tightness
\item and changes the contact structure.
\end{enumerate}

\smallskip
\noindent Proof of (3).

We will assume $\partial M=T$.  This is just one gluing theorem away from full generality.  The sutured manifold $(M,T)$, where $T$ is a toroidal suture, is automatically taut, and by Gabai's theorem it has a sutured manifold decomposition.  For simplicity, assume the first splitting surface $S$ intersects $T$ in a single curve, and say $(M,T) \stackrel{S}{\rightsquigarrow}(M', \gamma')$.

We would like to consider the corresponding convex decomposition, but first we must fix the toroidal suture.  Pick two parallel curves on $T$ dual to $S\cap T$ and define them to be $\Gamma$.  Figure~\ref{torusboundary} shows how we can force the usual correspondence between the sutured manifold decomposition on the first row and the convex structure on on the second row.  It is particularly important to note that $(S, \sigma)$ has boundary-parallel dividing curves.  Since $(M', \gamma')$ has a sutured manifold decomposition, $(M', \Gamma')$ carries a universally tight contact structure. By Theorem~\ref{Colin} $(M, \Gamma)$ does also.

The technique of adding a pair of parallel dividing curves to a boundary component with no sutures can be used in other settings as well.

\begin{figure}
    \centering
     \includegraphics{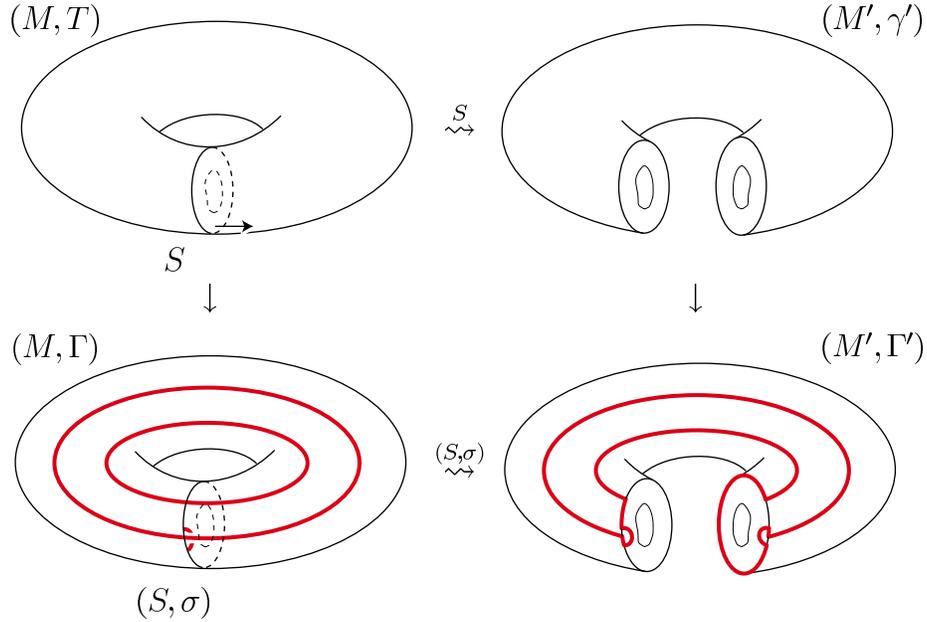}
     \caption{Introducing dividing curves on a torus suture.}
     \label{torusboundary}
\end{figure}

\smallskip
\noindent Proof of (4).

Continuing with the same $M, S$, and $T$, we need to show that the contact structure on $M \cup (T \times I)$ obtained by gluing the structure built in (3) and $\xi_k$ is universally tight.  Gluing along $T$ is beyond the scope of Theorem~\ref{Colin}.  Instead we will compare
$$
(M,\Gamma)\stackrel{(S,\sigma)}{\rightsquigarrow} (M\backslash S,\Gamma')
$$
and
$$
M\cup\xi_k\stackrel{S\cup A}{\rightsquigarrow} (M\backslash S) \cup (T\times I\backslash A)
$$
where $S$ is the first decomposing surface, and $A$ is an annulus extending $\partial S$ that is used to keep track of the $k$ twists in $\xi_k$.

The first row of Figure~\ref{Msplit} shows two views of $M\backslash S$ near $T\backslash S$.  In the first 3-dimensional picture,  a pair of dividing curves becomes a single dividing curve after corner rounding.  In the second,  the same neighborhood is expressed as a product with $S^1$, and the single dividing curve is shown as a point.  The second row of the figure gives a similar view of $(M\backslash S) \cup (T\times I\backslash A)$ in the case $k=1$.  The product with $S^1$ view also shows a convex surface transverse to the $S^1$ direction that detects the twisting along $\xi_k$.

\begin{figure}
    \centering
     \includegraphics{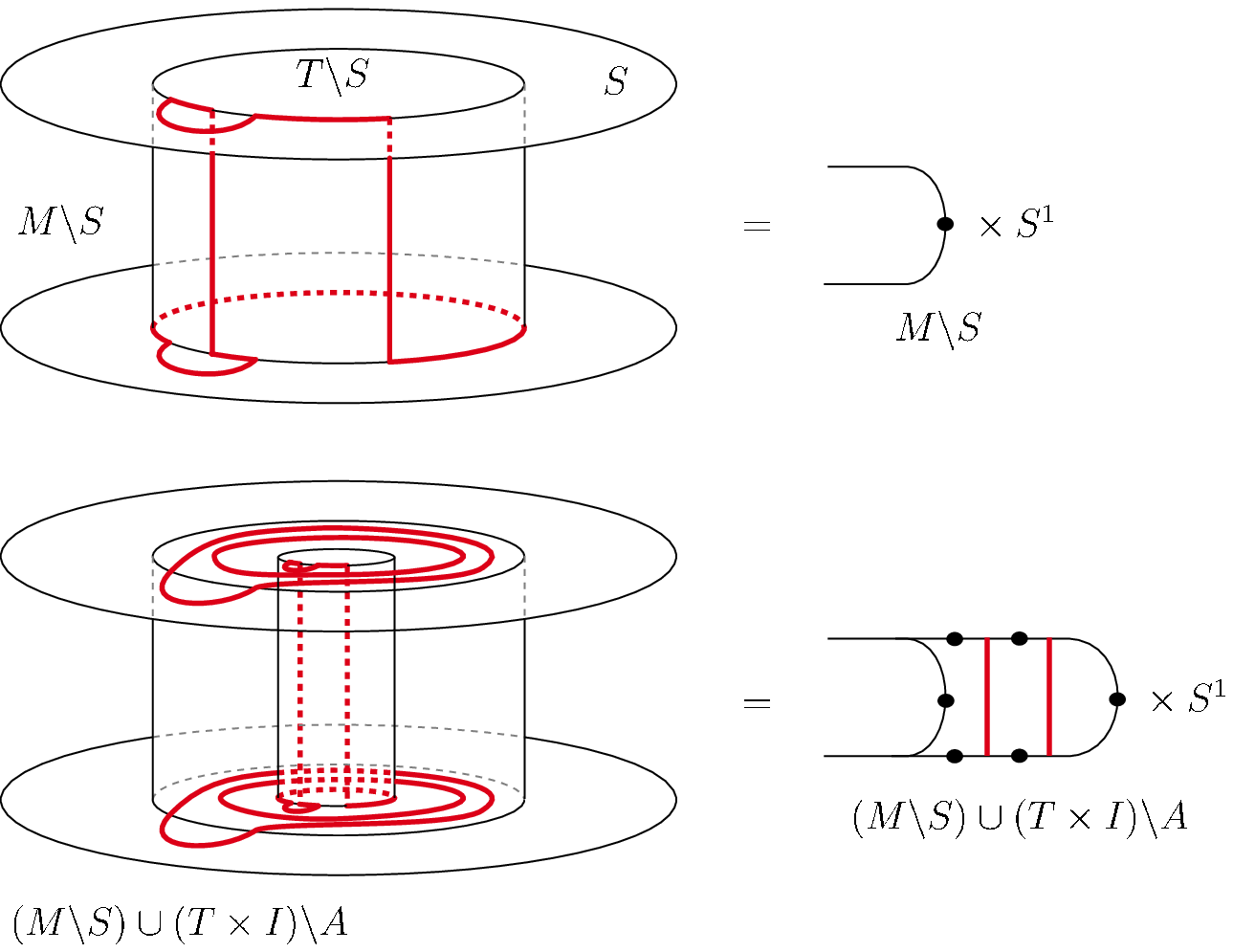}
     \caption{Views of $M\backslash S$ near the boundary before and after adding $T \times I \backslash A$.}
     \label{Msplit}
\end{figure}

We now show that our contact structure $(M\backslash S) \cup (T\times I\backslash A)$ is tight.  This will be done by finding an embedded copy of this space in the tight space $M\backslash S$, and then using the obvious (but useful!) fact that a subset of a tight space is tight.  Here is how this is done.

The curve $C$ parallel to $\partial S$ shown in Figure~\ref{isolating} is isolating.  Pass to a cover, without changing notation, in which the number of boundary components of $S$ is increased, and then $C$  becomes non-isolating.  Then use flexibility to make $C$ a Legendrian divide.  Figure~\ref{sideview} shows a product with $S^1$ view of $M\backslash S$ near $T\backslash S$.   The region shown is the product of a convex disk and $S^1$.  This disk is shown with part of two dividing curves that end on the two points of intersection with $C$, and these two points of $C$ are shown as hollow dots. 

\begin{figure}
    \centering
     \includegraphics{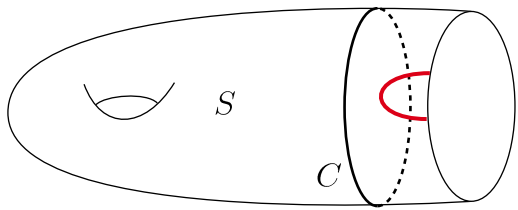}
     \caption{The curve $C$ is isolating.}
     \label{isolating}
\end{figure}

\begin{figure}
    \centering
     \includegraphics{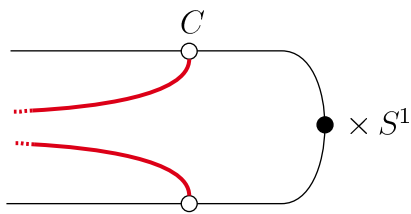}
     \caption{$C$ is now a Legendrian divide.}
     \label{sideview}
\end{figure}

Figure~\ref{Munionsideview} shows $(M\backslash S) \cup (T\times I\backslash A)$.  Finally, Figure~\ref{Bigsideview} shows a larger version of  $(M\backslash S)$ than Figure~\ref{sideview}.  The shaded subset of Figure~\ref{Bigsideview} is contact isomorphic to $(M\backslash S) \cup (T\times I\backslash A)$ and is necessarily (universally) tight.

\begin{figure}
    \centering
     \includegraphics{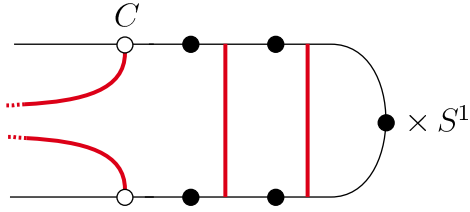}
     \caption{$M \cup (T \times I)$ split along $S \cup A$.}
     \label{Munionsideview}
\end{figure}

\begin{figure}
    \centering
     \includegraphics{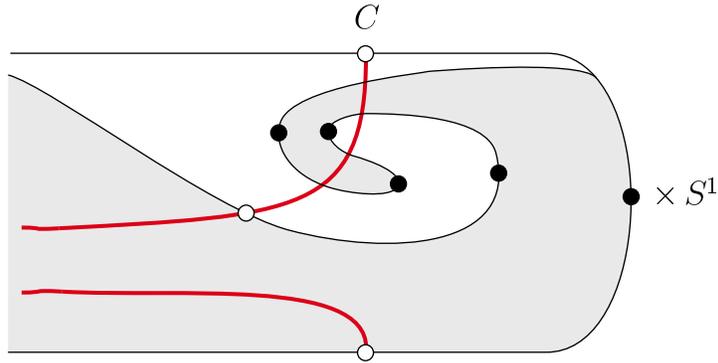}
     \caption{Enlarged view of $M \cup T$ split along $S$ shown with a distinguished subset.}
     \label{Bigsideview}
\end{figure}

Next we show the contact structure on $M \cup (T \times I)$ is tight.  This will require a gluing theorem along $S\cup A$, which unlike $S$, does not have boundary-parallel dividing curves.  Figure~\ref{SunionA} shows $\Gamma_{S \cup A}$.
\begin{figure}
    \centering
     \includegraphics{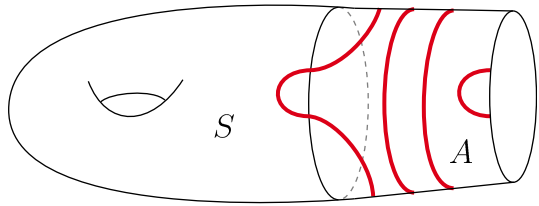}
     \caption{Dividing curves on $S \cup A$.}
     \label{SunionA}
\end{figure}

We will use the same gluing strategy:

\begin{itemize}
\item[(a)] Assume the union along $S\cup A$ is overtwisted.
\item[(b)] Isotop, via bypasses, $S \cup A$ off of the overtwisted disk, and
\item[(c)] argue that the split manifold stays tight during (b).
\end{itemize}

There is an important change in perspective though.  From the point of view of $S \cup A$, performing an isotopy  in (b) is equivalent to digging a bypass out of one side of $S \cup A$ and adding it to the other side.  We prove (c) by showing that there are no troublesome bypasses that can be dug from either side of $S \cup A$ in $M \cup (T \times I)$.

We have discussed bypasses which involve only boundary compressible dividing curves in the proof of Theorem~\ref{Colin}, so now we consider the existence of a bypass involving the closed dividing curves on $S \cup A$.

\begin{ex} Figure~\ref{bypasshere} shows a bypass attached along a dotted curve $\alpha$ connecting three different dividing curves on the boundary of $(M \cup T \times I)\backslash (S \cup A)$.   Figure~\ref{Munionsideview} showed the same space, but it showed a convex disk with a different set of dividing curves.  Certainly this bypass, if it exists, is not a subset of that convex disk.  Indeed the bypass itself may be very large and reach out of the portion of the manifold shown in either of these figures.
\end{ex}

\begin{figure}
    \centering
     \includegraphics{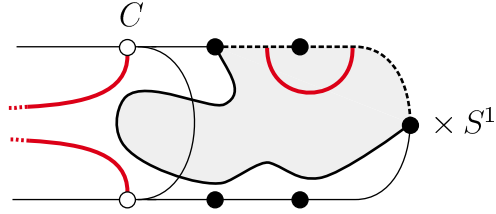}
     \caption{A possible location of a bypass.}
     \label{bypasshere}
\end{figure}

\begin{prop}
The bypass shown in Figure~\ref{bypasshere} can not exist.
\end{prop}

\begin{proof}
Consider the space obtained by adding the product of a bypass attached along $\alpha$ and $S^1$ that is shown in Figure~\ref{addbypass}.  The null-homotopic dividing curve that is created implies the contact structure is  overtwisted.

The opposite conclusion is reached if we add the same set and consider its relationship to the convex disk that we know is contained in $(M \cup T \times I)\backslash (S \cup A)$.  The union is shown in Figure~\ref{realeffect}.   We see that this space must be tight, for it too is a ``fold along $C$'', that is, it too can be found as a subset of $M\backslash S$.  This completes the proof of the Proposition and hence Step (4).
\end{proof}

\begin{figure}
    \centering
     \includegraphics{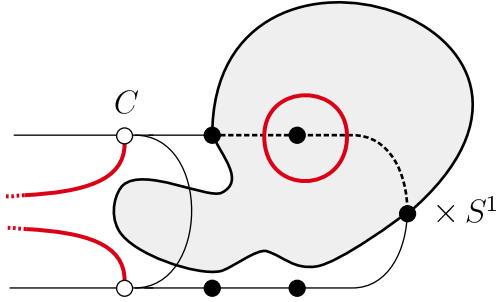}
     \caption{Detecting a bypass by adding a ``template''.}
     \label{addbypass}
\end{figure}

\begin{figure}
    \centering
     \includegraphics{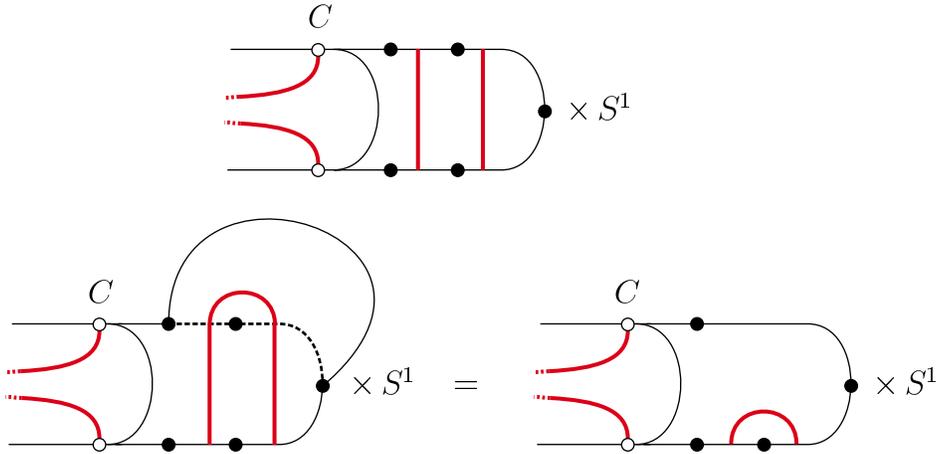}
     \caption{The actual result of adding the template.}
     \label{realeffect}
\end{figure}

\smallskip
\noindent Proof of (5).

We will only give the idea behind the technique used in producing invariants that distinguish the various contact structures on $M \cup (T \times I)$.  Let $F$ be a properly embedded convex surface which intersects the boundary component $T\times \{1\}$ of $M \cup (T \times I)$, and let $\delta$ be a homotopically essential arc in $F$ which  starts and ends on $T\times \{1\}$.  The minimum of $\#|\delta \cap \Gamma_F|$ over all such $\delta$ and $F$ is an invariant which tends to infinity as the twisting, that is, the $k$ in $\xi_k$, is increased.

\end{proof}

\section{Surface bundles}

We will give a classification of tight contact structures on $\Sigma \times I$ such that

\begin{itemize}

\item[(*)] $\Sigma$ is a closed surface with genus at least two, and the dividing curves on each component of $\partial(\Sigma \times I)$ are a pair of parallel non-separating curves. 

\end{itemize}

\begin{figure}
    \centering
     \includegraphics{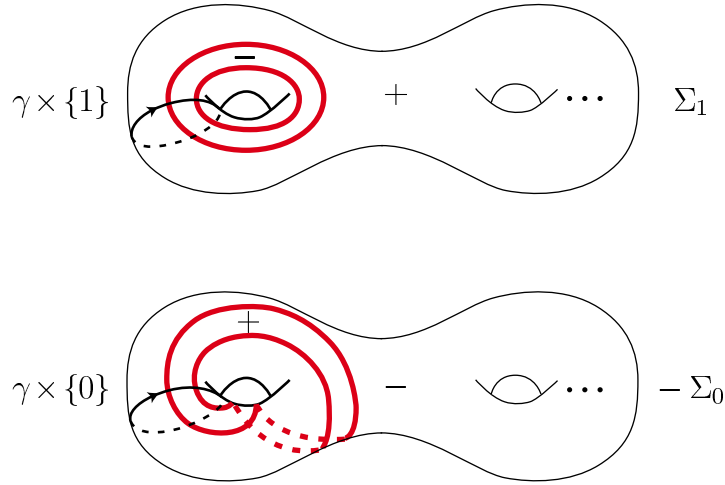}
     \caption{Dividing curves and splitting annulus on $\Sigma \times I$.}
     \label{fig2}
\end{figure}

\begin{ex}
In Figure~\ref{fig2}, notice that on $\Sigma_1$, $\chi(R_+)=\chi(\Sigma_1)$ and $\chi(R_-)=0$.  This is analogous to the product foliation on $\Sigma \times I$ in the sense that the contact 2-planes have outward pointing normal vectors everywhere on $\Sigma_1$, at least as measured by Euler charactertistic. This structure is a special case of an {\rm extremal} contact structure.
\end{ex}

\begin{defn}
A contact structure on a surface bundle with fibre $\Sigma$ is {\rm extremal}  if  the  Euler class of the 2-plane bundle, $e(\xi)$, satisfies $e(\xi)(\Sigma)=\pm\chi(\Sigma)$.  Equivalently, if $\Sigma$ is convex either $\chi (R_+) = \chi(\Sigma)$ or $\chi (R_-) = \chi(\Sigma)$.
\end{defn}

\begin{thm}\label{surface}\cite{HKM3}
There are exactly 4 (universally) tight non-product contact structures satisfying (*).  They correspond to a choice of dividing curve on each of $\Sigma_0$ and $\Sigma_1$.
\end{thm}

This theorem will be used to prove the following theorems.

\begin{thm}\label{bundles}
Let $\varphi$ be a pseudo-Anosov map of a closed surface $\Sigma$ . There is a unique extremal, tight (or universally tight) contact structure on $(\Sigma \times I)/(\varphi(x),1)\sim(x,0)$.
\end{thm}

\begin{thm}\label{GET}(Gabai-Eliashberg-Thurston Theorem)
If $M$ is Haken and $H_2(M)\ne 0$, then $M$ carries a universally tight contact structure.
\end{thm}

Theorem~\ref{GET} follows from Gabai's work \cite{Ga} on the existence of taut foliations  and Eliashberg and Thurston's perturbation technique \cite{ET} for producing universally tight contact structures from taut foliations.  The proof we give \cite{HKM4} is a direct construction which has the advantage of helping us discover new gluing theorems.

\begin{figure}
    \centering
     \includegraphics{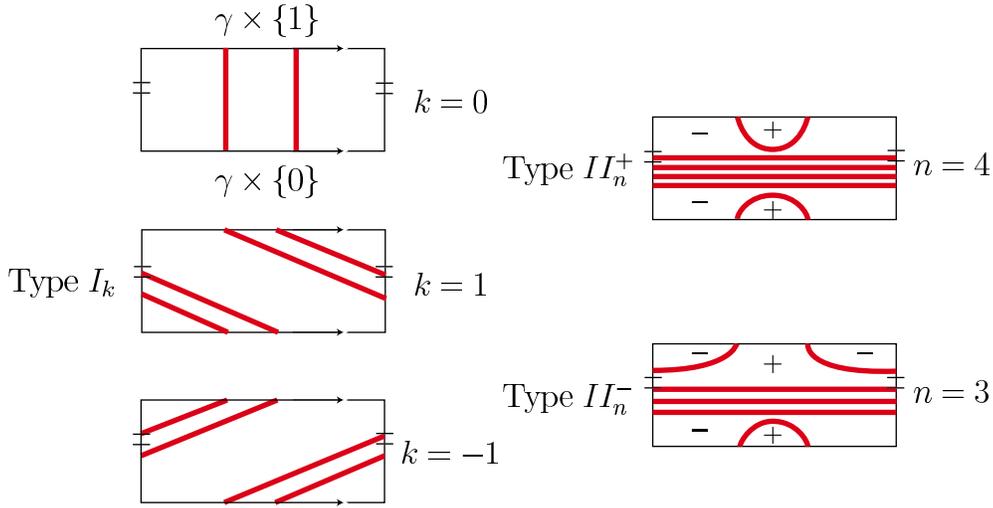}
     \caption{Possible dividing curves on $\gamma \times I$.}
     \label{fig3}
\end{figure}

\begin{proof}[Sketch of Theorem~\ref{surface}]
We concentrate only on the dividing curve configuration shown in Figure~\ref{fig2}.  What the proof strategy lacks in subtlety it makes up in directness.  We start by decomposing $\Sigma \times I$ along a vertical annulus $\gamma \times I$ whose boundary components are shown in Figure~\ref{fig3} and analyze all possible dividing curve configurations.

Figure~\ref{fig4} shows $\gamma \times I$ cut by a vertical arc into a rectangle, and it lists all possible dividing curve configurations such that the boundary components of $\gamma \times I$ intesect the dividing curves twice each.

\begin{figure}
    \centering
     \includegraphics{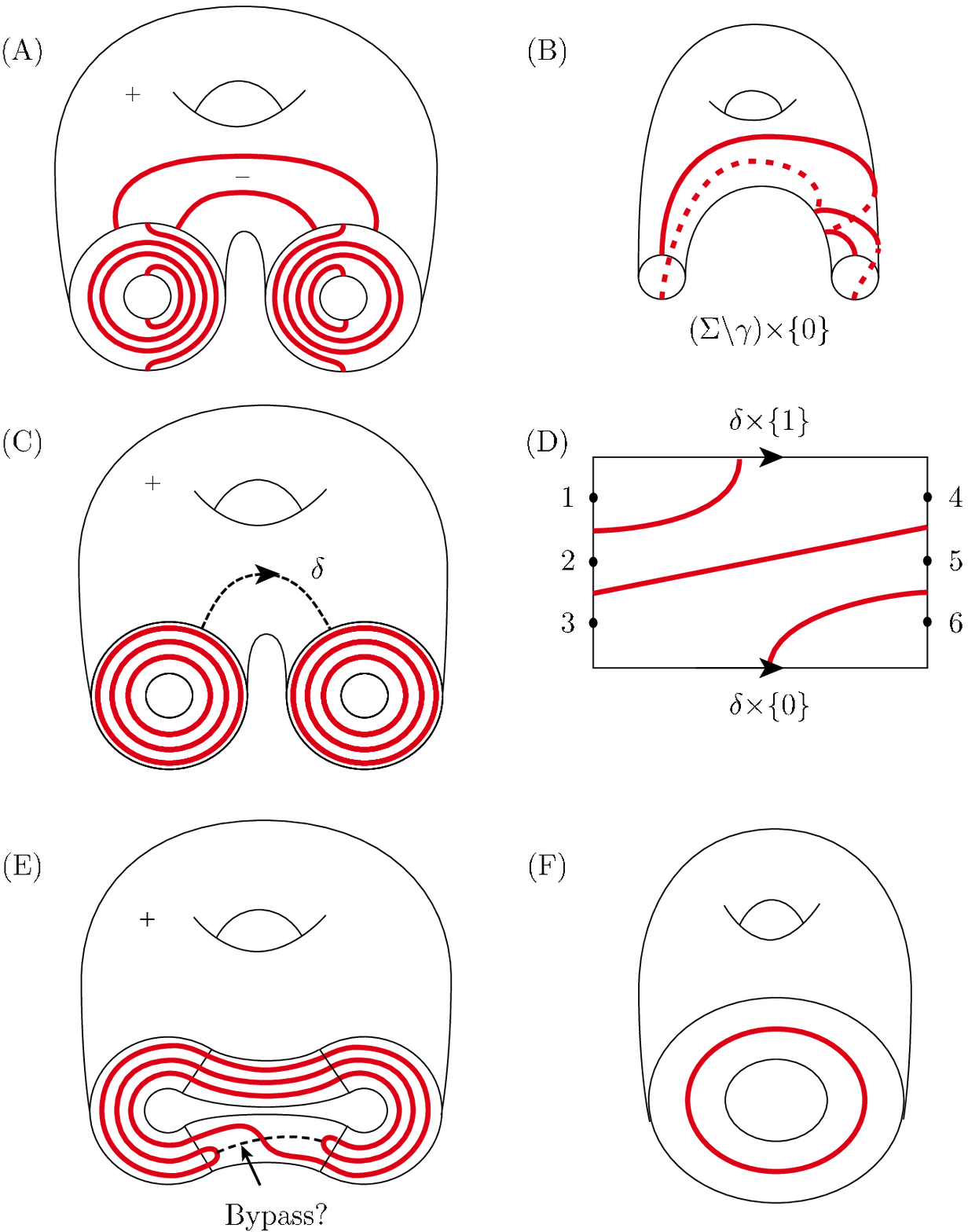}
     \caption{Decomposing $\Sigma \times I$.}
     \label{fig4}
\end{figure}

We will show how Type~$II_2^+$ can be reduced to Type~$II_0^+$, that is, if we start with a convex annulus of Type~$II_2^+$, we can find another convex annulus of Type~ $II_0^+$.  This case is fairly typical of the type of arguments we use to prove this classification theorem.

Figure~\ref{fig4}(A) shows the result of splitting $\Sigma \times I$ along an an annulus of Type~$II_2^+$, and Figure~\ref{fig4}(B) shows the dividing curves on $\Sigma \backslash \sigma \times \{0\}$.  After rounding the corners and gathering the dividing curves on the two vertical annuli of $(M \backslash \sigma) \times I$, the result is shown in Figure~\ref{fig4}(C).  Cut this along a convex rectangle $\delta \times I$, and again we must consider all possible dividing curve configurations on a splitting convex surface.

Suppose first that the dividing curve configuration on $\delta \times I$ has a boundary-parallel component whose half disk contains the point labelled 2 in $\delta \times I$. (Figure~\ref{fig4}(D) shows a different configuration.)  This would imply the existence of a bypass with arc of attachment running from 1 to 3.  Consider how this bypass is situated relative to $\sigma \times I$ -- it crosses two parallel dividing curves and ends on a third curve in Figure~\ref{fig4}(A).  But a quick computation (Example~\ref{parallel})  shows that the effect of pushing $\sigma \times I$ across such a bypass removes both closed dividing curves from $\sigma \times I$, that is, it produces a Type~$II_0^+$ annulus that we had promised to find.

The same logic shows that if any dividing curve is boundary-parallel and centered on the point 3, 4, or 5, a Type~$II_0^+$ can be shown to exist.  We are left with the dividing curve configuration of Figure~\ref{fig4}(D).  The result of cutting along $\delta \times I$ and rounding corners is shown in Figures~\ref{fig4}(E,F).

\begin{lemma}
There must exist a bypass along the arc of attachment shown in Figure~\ref{fig4}(E). \end{lemma}

\begin{proof}[Sketch of lemma]
If such a bypass exists exists then pushing the vertical annulus of Figure~\ref{fig4}(E) through it produces three dividing curves rather than just the one shown in Figure~\ref{fig4}(F).   One step in a complete proof is arguing that there is a vertical annulus with three dividing curves near the given vertical annulus.  This is very similar to the part of the  proof of tightness in the toroidal case, shown in Figure~\ref{Bigsideview}, in which a more complicated space was shown to exist inside a manifold with a single dividing curve on a vertical annulus in its boundary.  In short, the annulus is shown to exist by a folding argument along an isolating curve that is shown to exist near the original annulus.  This more complicated annulus implies the existence of the desired bypass.
\end{proof}

Pushing $\delta \times I$ through this bypass produces a new dividing curve configuration on $\delta \times I$ which has boundary-parallel dividing curves centered on the points labelled 3 and 4.  As above, this allows us to modify $\gamma \times I$ and produce an annulus of Type~$II_0^+$.

The rest of the proof of Theorem~\ref{surface} requires:
\begin{enumerate}

\item Many other similar reductions.

\item Existence of the four types of contact structures must be shown.

\item Uniqueness of these contact structures must be established.

\end{enumerate}

We  will show the existence and uniqueness of Type~$II_0^+$ tight contact contact structures.  Figure~\ref{typeII}(A, B) shows $\Sigma \times I$ cut along $\gamma \times I$.  Notice that the dividing curves on $\gamma \times I$ are boundary-parallel.  The result of rounding corners and the next splitting surface is shown in Figure~\ref{typeII}(C).  Since $\delta \times I$ intersects only two dividing curves, there is no choice for the diving curve set on $\delta \times I$, it must be the single boundary-parallel arc shown in Figure~\ref{typeII}(D).  The result of cutting along $\delta \times I$ and rounding corners is shown in Figure~\ref{typeII}(E,F).  And  now the pattern continues, $\epsilon \times I$ intersects dividing curves just twice, and continuing in this fashion we produce a convex decomposition in which all splitting surfaces have boundary-parallel dividing curves, and there is always a unique choice of tight dividing curves.

\begin{figure}
    \centering
     \includegraphics{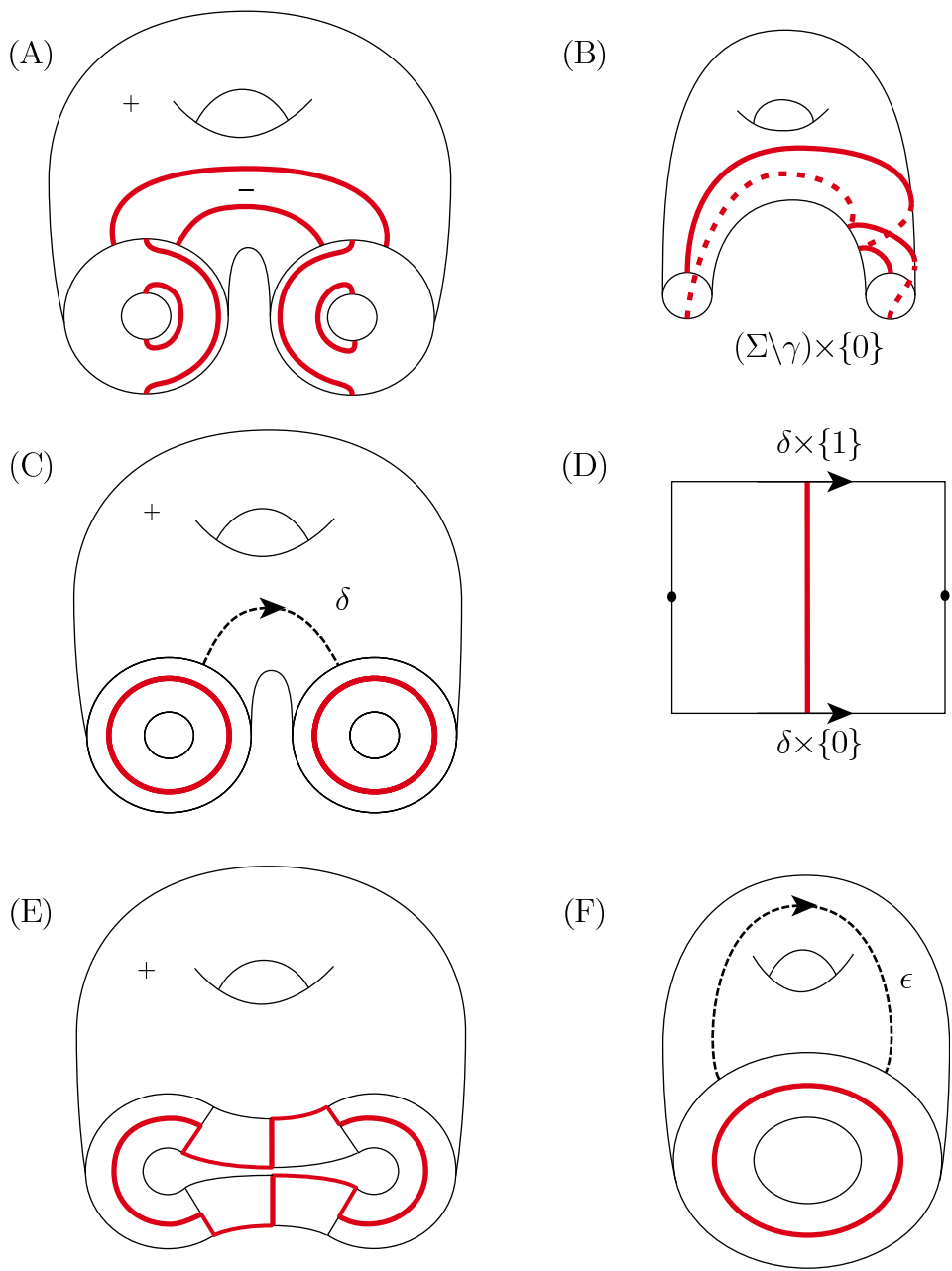}
     \caption{Existence of uniqueness of Type~$II^+_0$.}
     \label{typeII}
\end{figure}

Thus (2), in this case, follows from Theorem~\ref{Colin}, and (3) follows from the forced choice of dividing curves along the way.

It is worth emphasizing that a convex decomposition determines the contact structure near $\partial M$, near the splitting surfaces, and by uniqueness on all of the resulting $B^3$'s, that is, a convex decomposition determines the contact structure at every point of $M$.

This completes the sketch of Theorem~\ref{surface}.
\end{proof}

The statement of Theorem~\ref{surface} refers to a choice of dividing curves.  This choice can be described explicitly using the notion of a {\it straddled} dividing curve.

\begin{defn}
A dividing curve in $\partial (\Sigma \times I)$ is  {\rm straddled} if there exists a dual convex annulus with a boundary-parallel dividing curve centered on it.
\end{defn}

We record a couple of consequences of Theorem~\ref{surface} that will be used in applications.

\begin{cor}[Addition]
Let $\xi_1$ and $\xi_2$ be tight, non-product contact structures on $\Sigma \times [0, 1]$ and $\Sigma \times [1,2]$, respectively, that agree on $\Sigma \times \{1\}$.  Then $\xi_1\cup\xi_2$ is tight if and only if no dividing curve on $\Sigma \times \{1\}$ is straddled in both $\Sigma \times [0, 1]$ and $\Sigma \times [1,2]$. \qed
\end{cor}

\begin{cor}[Freedom of Choice]
Let $\xi$ be a non-product tight contact structure on $\Sigma \times I$, and let $\alpha_1$ and $\alpha_2$ be a pair of parallel, non-separating curves on $\Sigma$.  Then there exists a convex embedding of $\Sigma$ in $\Sigma \times I$ that is isotopic to the inclusion of a boundary component, such that $\Gamma_\Sigma = \alpha_1 \cup \alpha_2$. \qed
\end{cor}

The next proposition is self-evident and very useful \cite{H1}.

\begin{prop}[Imbalance Principle]
Let $S^1 \times [0,1]$ be a properly embedded convex annulus in $M$ such that $S^1 \times \{0\}$ intersects fewer dividing curves than $S^1 \times \{1\}$.  Then  $S^1 \times [0,1]$ contains a bypass centered on a dividing curve intersecing  $S^1 \times \{1\}$. \qed
\end{prop}

\begin{proof}[Sketch of Theorem~\ref{bundles}]

\begin{itemize}

\item Given a surface bundle with pseudo-Anosov holonomy $\varphi$, pick a fibre, isotop it until it is convex, and cut the bundle along the fibre.  The dividing curves on each boundary component consist of a family of parallel pairs of curves.

\item If there are more than one pair of parallel curves on either boundary component, then since $\varphi$ is pseudo-Anosov, there exists an imbalance annulus.

\item Isotoping the fibre through the bypass guaranteed by the Imbalance Principle reduces the number of dividing curves.  Continue until there is just a single pair on each boundary component.

\item By Freedom of Choice a new fibre can be chosen with a fixed pair of non-separating dividing curves.

\item Splitting the bundle along this fibre reduces an arbitrary bundle to one of the four standard forms given in Theorem~\ref{surface}.  Of the four possible straddlings, two are ruled out because of the tightness of the gluing that recreates the original surface bundle.  The other two are related by another application of Freedom of  Choice.

\end{itemize}

\end{proof}

\begin{proof}[Sketch of Theorem~\ref{GET}]

By Theorem~\ref{equivalent} we may assume $M$ is a closed manifold.  Let  $\Sigma\subset M$ be a Thurston norm minimizing surface corresponding to a non-zero element of $H_2(M)$ and split  $M$ along $\Sigma$.

\begin{figure}
    \centering
     \includegraphics{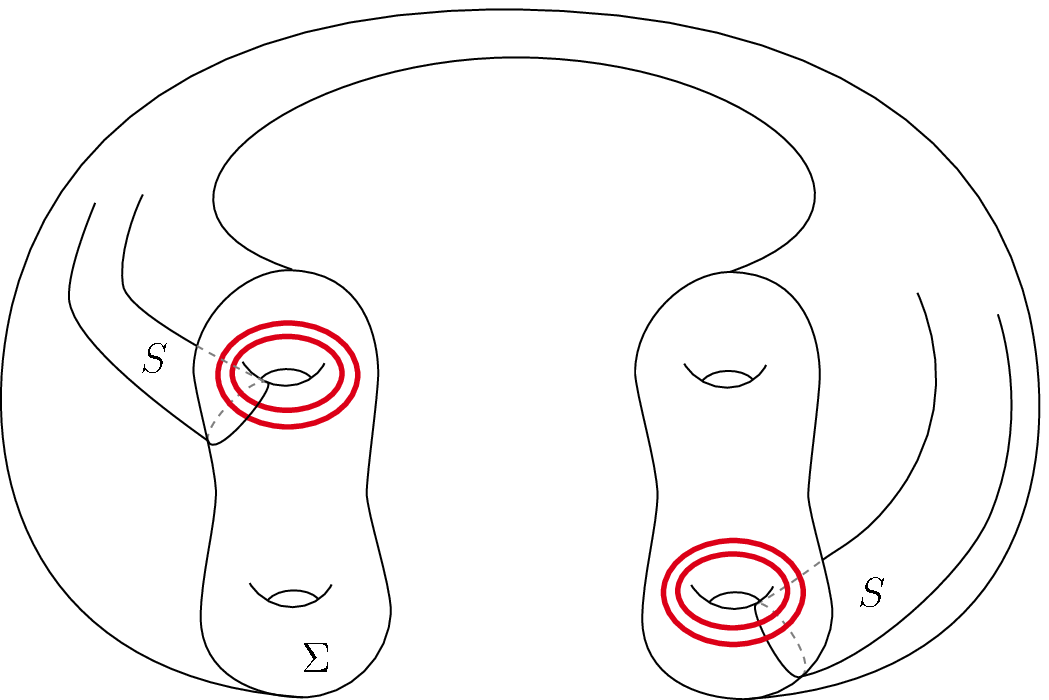}
     \caption{The first splitting surface $S$ in $M\backslash\Sigma$.}
     \label{GETdecomp}
\end{figure}

The sutured manifold $M\backslash\Sigma$ has no sutures, but it does have a sutured manifold decomposition.   Let the first splitting surface be $S$, and we shall assume that $S$ intersects each copy of $\Sigma$ in a single closed curve as shown in Figure~\ref{GETdecomp}.  Make $M \backslash \Sigma$ a convex structure by adding a pair of parallel dividing curves dual to $\partial S$ on each boundary component.  Make $S$ a convex surface by adding boundary compressible dividing curves $\sigma$ straddling a component of $\Gamma$, the dividing set on the boundary of $M\backslash \Sigma$, on each copy of $\Sigma$.  If the right curves are straddled, splitting the sutured manifold $M \backslash \Sigma$ along $S$ will correspond to the convex splitting defined by $(S, \sigma)$.  The remaining steps of the convex decomposition are directly inherited from the sutured manifold decomposition.

This convex decomposition of $(M\backslash\Sigma, \Gamma)$ is by surfaces all of whose dividing curves are boundary-parallel.  Thus by Theorem~\ref{Colin}, there is a tight contact structure on $M\backslash\Sigma$.  Moreover, by construction, one dividing curve on each boundary component of $M \backslash \Sigma$ is straddled by a dividing curve on $S$.

By Theorem~\ref{surface} there are four choices of tight contact structure on $\Sigma \times I$ that could be used to attach to $M \backslash \Sigma$ and produce a tight contact structure on $M$.  It should seem very plausible, and it is true, that a curve straddled on both sides gives rise to an overtwisted disk.  Thus we insert the unique, non-product, contact structure on $\Sigma \times I$ that gives $(M \backslash\Sigma) \cup\Sigma \times I$ a chance of being tight.

\begin{figure}
    \centering
     \includegraphics{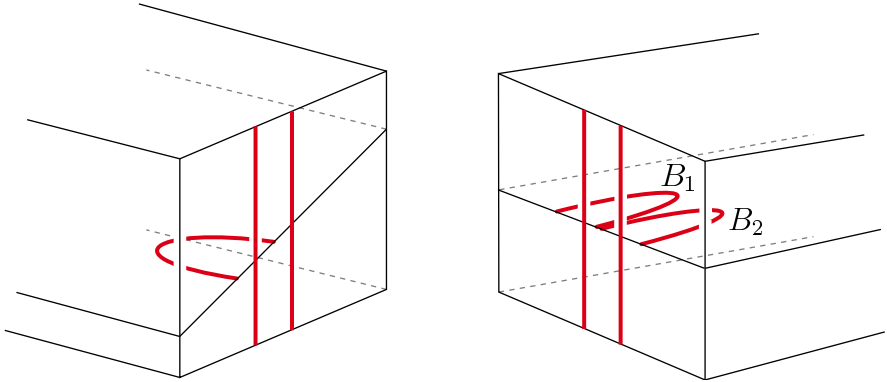}
     \caption{Possible location of bypasses before gluing.}
     \label{glue}
\end{figure}

We are two gluing theorems away from a complete proof of tightness on $M$; we must glue along each of the boundary components of $\Sigma \times I$.   As we have seen, the general form of these gluing proofs is: 

\begin{enumerate}

\item Given an overtwisted disk in $M$, push $\Sigma$ off it using bypasses while keeping $M\backslash\Sigma$ tight.

\item Analyze which bypasses exist on one component of $M \backslash \Sigma$ and which can be added to the other component while preserving tightness.

\end{enumerate}

Rather than do this in generality, consider the local version of this that is shown in Figure~\ref{glue}.  On the left are two dividing curves one of which is straddled.  On the right are the two dividing curves about to be identified with the curves on the left.  Also shown are two bypasses.  The first, $B_1$, is known to exist by construction, thus if it is removed and added to the other side tightness must be shown to be preserved.  The second, $B_2$, if added to the left would produce an overtwisted disk, thus, as part of a sufficient gluing theorem, these must be shown not to exist.  These local gluing results follow from the next lemma.

\end{proof}

\begin{lemma}
Let $\gamma_u$ and $\gamma_s$ be a pair of parallel dividing curves on $\partial M$, and assume $\gamma_s$ is straddled and the contact structure on $M$ is tight.  Then adding a bypass to $M$ across $\gamma_u$ will produce a tight contact structure.
\end{lemma}

Since adding a bypass to $M$ across $\gamma_s$ produces an overtwisted structure, it follows that $\gamma_u$ is not straddled.

\begin{proof}
Figure~\ref{bypasses} shows a neighborhood $A \times I$ of an  annular neighborhood $A$ of $\gamma_s$ and $\gamma_u$ in $\partial M$.  It also shows the arc of attachment $\alpha$ which straddles $\gamma_s$ and the arc of attachment $\beta$ to which a bypass is being added.  The annulus parallel and below $A$ shows the result of removing the bypass attached along $\alpha$, and the annulus above $A$ shows the result of adding a bypass along $\beta$.  The figure on the right shows the dividing curves on the boundary $A \times I$. 

At least on the boundary, the figure on the right looks like a product contact structure on $A \times I$, and indeed it is.  Since the attaching curves, $\alpha$ and $\beta$, are disjoint, the contact structure on $A$ can be built by first attaching a bypass to the bottom annulus along $\beta$ and then attaching a bypass along $\alpha$.  From this point of view, the isotopy class of the dividing curves remains unchanged after adding each bypass (see also Example~\ref{trivial}), and thus the contact structure is a product.  It now follows that adding a bypass across $\beta$ is the same as removing a bypass in $M$ attached along $\alpha$, and this operation preserves tightness.

\begin{figure}
    \centering
     \includegraphics{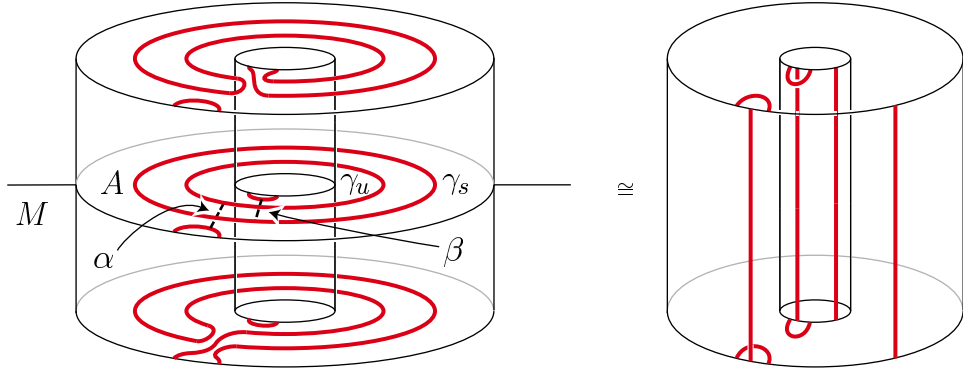}
     \caption{Adding a bypass across $\beta$ is the same as removing a bypass across $\alpha$.}
     \label{bypasses}
\end{figure}

\end{proof}

\section{Open Questions}

There are two fundamental classes of open questions:

\begin{enumerate}

\item Which $M^3$ carry tight contact structures?

\item What are the topological implications of carrying a tight contact structure?

\end{enumerate}

The central existence question, particularly from the point of view developed in this paper, is the question of whether or not Haken homology spheres $M$ always carry tight contact structures.

\begin{figure}
    \centering
     \includegraphics{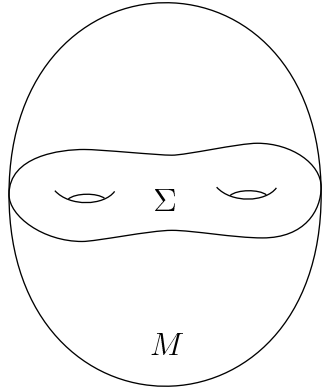}
     \caption{Haken homology sphere.}
     \label{haken}
\end{figure}

In such a manifold, every surface $\Sigma\subset M$ must separate.  In particular if there is a tight contact structure $\xi$ on $M$, then $e(\xi)(\Sigma)=0$.  This means that if $\Sigma$ is convex, then $\chi(R_+) = \chi(R_-)$.  This is exactly opposite to the extremal case when $\chi(R_+) = \chi(\Sigma)$ and $\chi(R_-)=0$.

Presumably constructing contact structures will involve:

\begin{itemize}

\item Classification of such structures on $\Sigma \times I$ and

\item new gluing theorems.

\end{itemize}

\begin{figure}
    \centering
     \includegraphics{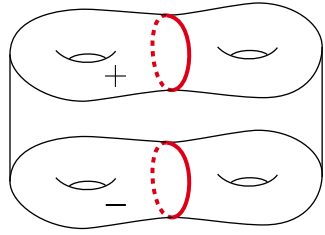}
     \caption{A non-extremal boundary configuration.}
     \label{onedividing}
\end{figure}

\begin{ex}
Perhaps the simplest example of this sort of classification question on $\Sigma \times I$ is shown in Figure~\ref{onedividing}.  Preliminary work of Cofer \cite{Co} shows there is exactly one tight, non-product, contact structure with these dividing curves.  This example has the bizarre property that if you add any non-trivial bypass, it becomes overtwisted.  It follows that it does not occur as a subset of any tight contact structure on $\Sigma \times I$ other than itself, and it may not show up in any tight closed 3-manifold.
\end{ex}

Very little is known about (2), implications of carrying a tight contact structure, so we will describe results that have been obtained in lamination theory, that perhaps have analogues in contact topology.

\begin{defn}
A {\rm lamination} of $M^3$ is a disjoint union of surfaces which are locally homeomorphic to the product of $D^2$ and a closed subset of $I$.  

A lamination is {\rm essential} if the leaves are incompressible, the complementary regions are irreducible, and there are no folded leaves. 

A lamination is {\rm genuine} if it is essential and some complementary region is not a product of a boundary leaf and $I$.
\end{defn}

\begin{figure}
    \centering
     \includegraphics{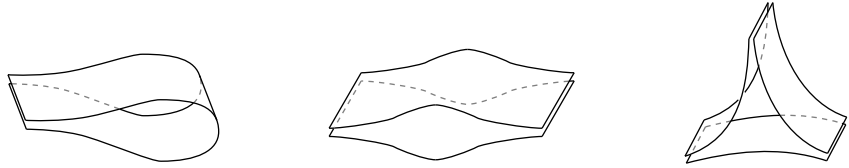}
     \caption{Complementary regions of a general lamination.}
     \label{essential}
\end{figure}

Figure~\ref{essential} shows, in order, a folded leaf, a complementary region that is a product of a boundary leaf and $I$, and a complementary region that is not such a product. 

\begin{defn}
The {\rm Euler characteristic} of a surface with cusped boundary is defined to be the usual Euler characteristic of the underlying space minus half of the number of cusps.
\end{defn}

The cross-sections of the complementary regions shown in Figure~\ref{essential} are a disk with one cusp ($\chi = 1/2)$, a disk with two cusps ($\chi = 0$), and a disk with three cusps ($\chi = -1/2$).  The definition of essential consists of bans on various types of positive Euler characteristic, while the notion of a genuine lamination postulates the existence of some negative Euler characteristic in $M$.  We shall see that atoroidal manifolds are group negatively curved (Theorem~\ref{gnc}).  It is not clear what additional structure should be made for contact structures that might make the this theorem apply in that setting as well.

\begin{figure}
    \centering
     \includegraphics{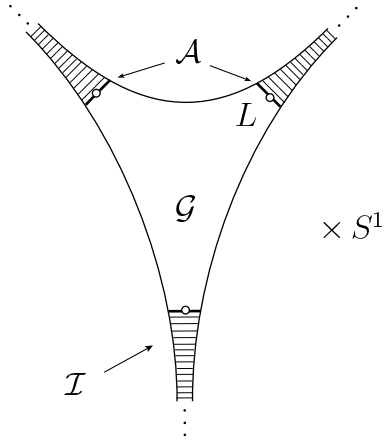}
     \caption{Structure on the complementary region of a genuine lamination.}
     \label{genuine}
\end{figure}

By the JSJ decomposition theorem, there is a unique $I-$bundle structure $\mathcal I$ on the ends of each complementary region.  Thus each complementary region decomposes as the union of a $\mathcal I$ and the {\it guts} $\mathcal G$ as shown in Figure~\ref{genuine}.

The key features of this decomposition are:

\begin{itemize}
\item $\mathcal G$ is compact.

\item By maximality of $\mathcal I$, $\mathcal G$ has no {\it product disks}, that is, there are no non-trivial rectangles in $\mathcal G$ with sides that alternately consist of $I$-bundle fibres of $\mathcal I$ and arcs in leaves of the lamination.

\item An essential lamination is genuine if and only if ${\mathcal G}\ne\emptyset$.

\item ${\mathcal G}\cap{\mathcal I}$ is a finite union of annuli $\mathcal A$.  The union of the cores of these annuli are a link denoted $L$.
\end{itemize}

\begin{defn}  
$M$ is {\rm group negatively curved} if there exists a constant $C$ such that for every null-homotopic curve, $f:S^1 \to M$, there exists an extension of $f$ to a disk $D$ such that 
\[
	\mbox{area}(f(D)) < C \cdot \mbox{length} (f(\partial D)).
\]

$M$ is {\rm group negatively curved with respect to a link $L$ in $M$} if there exists a constant $C$ such that for every null-homotopic curve $f:S^1 \to M$, there exists an extension of $f$ to a disk $D$ such that 
\[
	\mbox{area}(f(D)) < C \cdot (\mbox{length} (f(\partial D)) + \mbox{wr} (f(\partial D), L)).
\]

The {\rm wrapping number} $\mbox{wr}(f(\partial D), L)$ is a geometric linking number and is defined to be the minimum, taken over all disks $E$ with $\partial E = f(\partial D)$ of the number of points of intersection of $E$ with $L$.
\end{defn}

The inequality in the definition of {\it group negatively curved with respect to a link $L$ in $M$} is equivalent to the existence of a constant such that at least one of the two inequalities is satisfied:
\[
	\mbox{area}(f(D)) < 2C \cdot \mbox{length} (f(\partial D))
\]
or
\[
	\mbox{area}(f(D)) < 2C \cdot  \mbox{wr} (f(\partial D), L).
\]

We will need the following remarkable theorem.

\begin{thm}\label{gut}(Gabai's Ubiquity Theorem \cite{Ga98})
If $M$ is closed, irreducible, and atoroidal, and if $L \not \subset B^3$, then $M$ is group negatively curved with respect to  $L$. \qed
\end{thm}

\begin{thm}\label{gnc}\cite{GK2}
If $M$ is atoroidal and contains a genuine lamination $\lambda$, then $M$ is group negatively curved.
\end{thm}

Before applying Gabai's Ubiquity Theorem to the proof of this, we will need the following lemma which says that to prove an isoperimetric inequality  for all null-homotopic curves, it is enough to prove the inequality on a ``dense'' subset.

\begin{lemma}
Let $A$ be the set of all null-homotopic curves $g:S^1 \to M$, and let $S$ be a subset of $A$.  If

\begin{itemize}
\item all $f\in S$ satisfy an isoperimetric inequality,

\item each $g\in A$ is approximated by an $f\in S$ by a small area homotopy, and

\item $\mbox{length}(f)$ is not drastically longer than $\mbox{length}(g)$,

\end{itemize}
then all $g\in A$ satisfy an isoperimetric inequality.
\end{lemma}

\begin{proof} This follows by piecing together the homotopies shown in Figure~\ref{small}.

\begin{figure}
    \centering
     \includegraphics{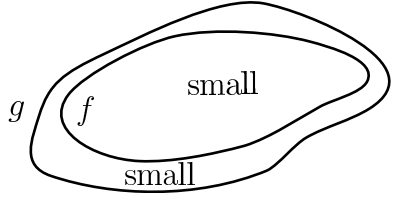}
     \caption{A small area homotopy that doesn't increase length much.}
     \label{small}
\end{figure}
\end{proof}

\begin{proof}[Sketch of Theorem~\ref{gnc}]
To apply this lemma think of $\mathcal G$ as a big, fat subset of $M$.  Then to show $M$ is group negatively curved, it is enough to prove an isoperimetric inequality for the set of null-homotopic curves $f:S^1 \rightarrow M$ such that

\begin{enumerate}
\item $f$ is transverse to  $\lambda$.

\item Each component of $f^{-1}({\mathcal G})$ has length greater than some constant $\varepsilon$.
\end{enumerate}
In other words, short bits of $f^{-1}({\mathcal G})$ can be efficiently removed as in Figure~\ref{shortbits}.

\begin{figure}
    \centering
     \includegraphics{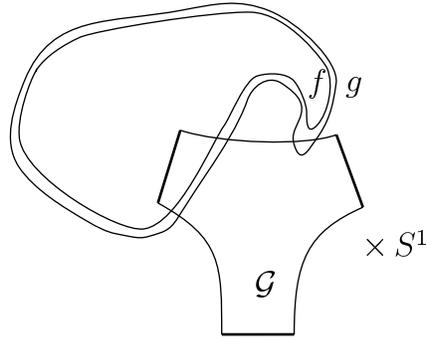}
     \caption{Short portions of $g^{-1}(\mathcal G)$ can be removed efficently.}
     \label{shortbits}
\end{figure}

Given such a null-homotopic $f:S^1 \to M$, there is, by Theorem~\ref{gut}, a disk of null-homotopy, $D$, such that 
\[
	\mbox{area}(f(D)) < 2C \cdot | f(D) \cap L |. 
\]
Thus it will be enough to show that for some constant $C'$
\[
	2C \cdot | f(D) \cap L | < C' \cdot \mbox{length}(\partial f(D)).
\]

Figure~\ref{pullback} shows $f^{-1}(D)$.  The figure shows $f^{-1}(\mathcal G)$ as shaded, and $f^{-1}(\mathcal I)$ as white.  Since we are only trying to give a sketch of the main ideas, we will think of $f$ as an embedding.

\begin{figure}
    \centering
     \includegraphics{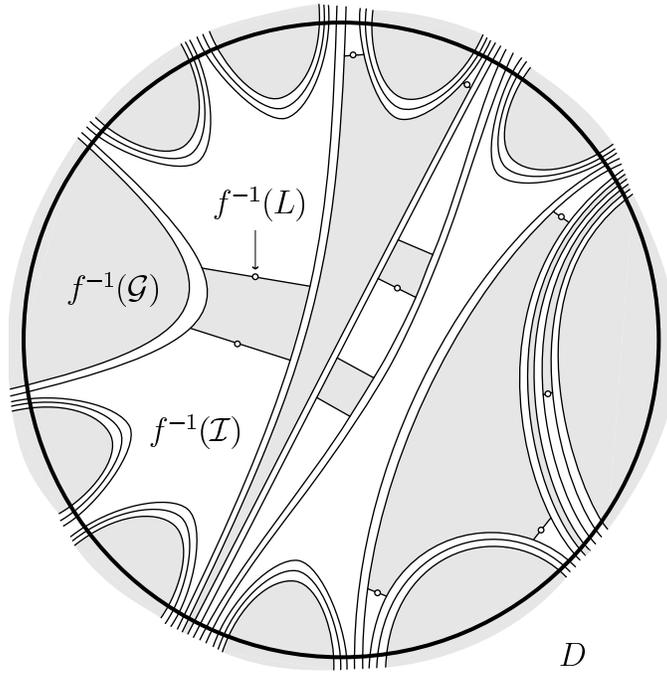}
     \caption{The pullback of $\lambda,  {\mathcal G}$, and ${\mathcal I}$ to $D$.}
     \label{pullback}
\end{figure}

Figure~\ref{positive} shows regions that might occur as subsets of $f^{-1}(D)$.  The first region, a null-homotopic circle, can be removed by choosing a new map of $D \to M$ since leaves of $\lambda$ are incompressible.  The second region, a folded leaf, can not occur in an essential lamination.  And finally the third region, a half disk mapped into $\mathcal G$ does occur, and thus we arrive at 

\smallskip
\noindent {\it Conclusion 1.} Regions of  $f^{-1}(D)$ with positive Euler characteristic contribute at least $\varepsilon$ to the length $f(\partial D)$.
\smallskip

\begin{figure}
    \centering
     \includegraphics{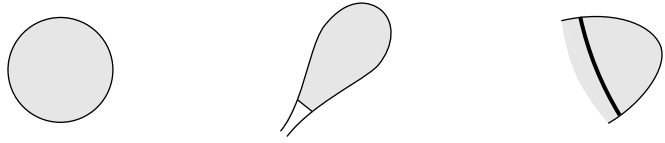}
     \caption{Possible regions of $D$ with positive Euler characteristic.}
     \label{positive}
\end{figure}

Figure~\ref{nonpositive} shows typical regions of $f^{-1}(D)$ which contain points of $f^{-1}(L)$.  The first figure, a cusped triangle, has negative Euler characteristic.  The second region shown has Euler characteristic zero and is a product disk in $\mathcal G$.  This can not exist by the definition of $\mathcal G$.  The third region also has Euler characteristic zero, but it contains an arc that is mapped into $\mathcal G$, thus it contributes at least $\varepsilon$ to length$(f(\partial D))$.  After removing the middle regions that can not exist we reach

\begin{figure}
    \centering
     \includegraphics{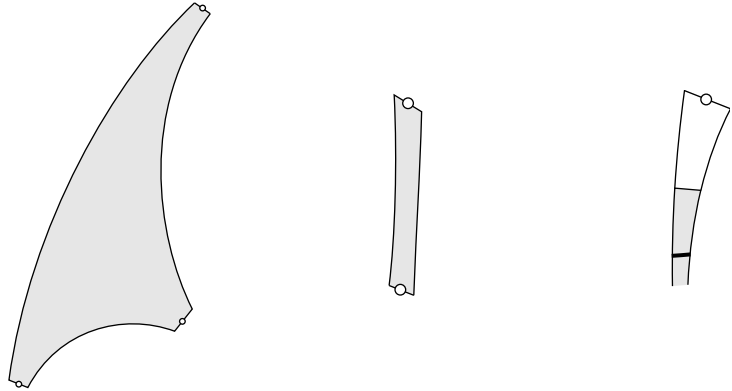}
     \caption{Regions of $D$ which contain points of $f^{-1}(L)$.}
     \label{nonpositive}
\end{figure}

\smallskip
\noindent {\it Conclusion 2.}  Points of $f^{-1}(L)$ either show up in regions of negative negative Euler chacateristic or  they contribute at least $\varepsilon$ to length$(f(\partial D))$.
\smallskip

We can now complete the proof.  We have a disk $D$ such that
\[
	\mbox{area}(f(D))<2C \cdot | f(D) \cap L |,
\]
thus a large area disk gives many points of $f^{-1}(L)$.  By Conclusion~2, these points either directly contribute to the length of $f(\partial D)$, or they show up in regions of negative Euler characteristic.  But $\chi(D)=1$, thus the existence of regions with negative Euler characteristic implies the existence of regions with positive Euler characteristic.  By Conclusion~1,  these in turn contribute even more to the length of $f(\partial D)$.  Thus we conclude
\[
	\mbox{area}(f(D))<2C \cdot | f(D) \cap L | < C' \cdot \mbox{length}(f(\partial D)).
\]

\end{proof}

A key feature of this proof that doesn't have an obvious analogue in contact topology is the crude notion of length given by pulling back $\mathcal G$ to $\partial D$.

We would like to end up by pointing out that there are not clear connections between tight contact structures on $M$ and the fundamental group of $M$.  For instance, it is not known if a homotopy 3-sphere supports a tight contact structure whether it must be $S^3$.

By way of contrast, there are many $\pi_1(M)$ actions that can be constructed from foliations and laminations.  The {\it leaf space} is the quotient of the universal cover by leaves and complementary regions.  The quotient is an {\it order tree}, and there is always an action of $\pi_1(M)$ on it.

Bestvina and Mess \cite{BM} show that if $M$ is group negatively curved then there is an action of $\pi_1(M)$ on $S^2$.  This can be applied to the manifolds of Theorem~\ref{gnc}, and indeed by Calegari's work \cite{Ca}, there are far more manifolds in this collection than originally realized.

Palmeira's Theorem \cite{Pa} is generalized to laminations in \cite{GK1}, and it follows that the universal cover $(\tilde M, \tilde{\lambda})$ is always homeomorphic to a product $({\mathbb R}^2, \kappa) \times {\mathbb R}$ where $\kappa$ is a lamination of the plane.  Calegari and Dunfield \cite{CD} point out that $({\mathbb R}^2, \kappa)$ can be thought of as $({\mathbb H}^2, \kappa)$ and from this they can sometimes produce an action on $S_{\infty}^{1}$.

Calegari and Dunfield \cite{CD} have more general results.  They generalize Thurston's work on the universal circle, and using Candel's theorem \cite{Can}, they identify leaves of $\lambda$ with ${\mathbb H}^2$, and they identify all $S_{\infty}^{1}$'s coming from the ${\mathbb H}^2$'s to get a $\pi_1(M)$ action on $S^1_{\smaller\smaller\mbox{univ}}$.  This works for taut foliations and some genuine laminations.

\smallskip

{\it Acknowledgements.}   This article is based on a series of talks given at the Tokyo Institue of Technology from June~3--7, 2002.  I would like to thank Professor Matsumoto for inviting me  and Professor Kojima for hosting my visit.  I would particularly like to thank Ko Honda for his lectures during an informal seminar at the University of Georgia that introduced me to his point of view in contact topology and for suggesting many improvements to this paper.

\end{document}